\documentclass[11pt]{article}
\usepackage{tocbibind} % adds References to the pdf bookmarks, in amsart this seems to be standard, but not in article
\usepackage{amsmath,amsthm,amssymb,thmtools}
\usepackage{fancyhdr}
\usepackage[british]{babel}
\usepackage{geometry,mathtools}
\usepackage{enumitem}
\usepackage{algpseudocode}
\usepackage{dsfont}
\usepackage{centernot}
\usepackage{xstring}
\usepackage{colortbl}
\usepackage{pifont}

\usepackage[section]{algorithm}
\usepackage{graphicx}
\usepackage[font={footnotesize}]{caption}
\usepackage[usenames,dvipsnames,table]{xcolor}
\usepackage{pstricks,tikz}
\usetikzlibrary{patterns,calc,fadings,decorations.pathreplacing}
\usepackage[h]{esvect}
\usepackage[
bookmarksopen=true,
bookmarksopenlevel=1,
colorlinks=true,
linkcolor=darkblue,
linktoc=page,
citecolor=darkblue,
urlcolor=darkblue
]{hyperref}
\usepackage{bbm}

\usepackage{multirow}

\numberwithin{equation}{section}

\definecolor{darkblue}{rgb}{0,0,0.5}

\newdimen\margin
\def\textno#1&#2\par{
   \margin=\hsize
   \advance\margin by -4\parindent
          \setbox1=\hbox{\sl#1}
   \ifdim\wd1 < \margin
      $$\box1\eqno#2$$
   \else
      \bigbreak
      \hbox to \hsize{\indent$\vcenter{\advance\hsize by -3\parindent
      \it\noindent#1}\hfil#2$}
      \bigbreak
   \fi}

\usepackage{cleveref}[capitalise, compress, nameinlink, noabbrev]

\theoremstyle{plain}
\newtheorem{theorem}{Theorem}[section]
\crefname{theorem}{Theorem}{Theorems}

\newtheorem{prop}[theorem]{Proposition}
\crefname{prop}{Proposition}{Propositions}

\crefname{cor}{Corollary}{Corollaries}

\newtheorem{cl}[theorem]{Claim}
\crefname{cl}{Claim}{Claims}

\newtheorem{lemma}[theorem]{Lemma}
\crefname{lemma}{Lemma}{Lemmata}
\newtheorem{conj}[theorem]{Conjecture}
\crefname{conj}{Conjecture}{Conjectures}

\crefname{problem}{Problem}{Problems}

\crefname{observation}{Observation}{Observations}

\crefname{setup}{Setup}{Setups}

\crefname{fact}{Fact}{Facts}

\crefname{remark}{Remark}{Remarks}

\crefname{example}{Example}{Examples}

\theoremstyle{definition}
\newtheorem{defin}[theorem]{Definition}
\crefname{defin}{Definition}{Definitions}

\crefname{construction}{Construction}{Constructions}

\crefname{question}{Question}{Questions}

\numberwithin{equation}{section}

\crefname{section}{Section}{Sections}
\crefname{appendix}{Appendix}{Appendix}

\crefname{figure}{Figure}{Figures}

\def\proof{\removelastskip\penalty55\medskip\noindent\begin{stepenv}\end{stepenv}{\bf Proof. }} % in each main proof, claim and step counter set back

\def\lateproof#1{\removelastskip\penalty55\medskip\noindent\begin{stepenv}\end{stepenv}{\bf Proof of #1. }} % in each main proof, claim and step counter set back

\def\noproof{{\unskip\nobreak\hfill\penalty50\hskip2em\hbox{}\nobreak\hfill%
       $\square$\parfillskip=0pt\finalhyphendemerits=0\par}\goodbreak}
\def\endproof{\noproof\bigskip}

\newcounter{stepenv}
\newenvironment{stepenv}[1][]{\refstepcounter{stepenv}}{}

\newcounter{step}[stepenv]

\newcounter{substep}[step]
\renewcommand{\thesubstep}{\thestep.\arabic{substep}}

\newcounter{claim}[stepenv]

\newcommand{\cA}{\mathcal{A}}
\newcommand{\cB}{\mathcal{B}}
\newcommand{\cC}{\mathcal{C}}
\newcommand{\cD}{\mathcal{D}}
\newcommand{\cE}{\mathcal{E}}
\newcommand{\cF}{\mathcal{F}}

\newcommand{\cL}{\mathcal{L}}

\newcommand{\cP}{\mathcal{P}}

\newcommand{\cS}{\mathcal{S}}

\newcommand{\bR}{\mathbb{R}}

\def\eps{{\varepsilon}}
\def\e{{\mathrm{e}}}

\newcommand{\defn}{\emph}

\newcommand{\prob}[1]{\mathrm{\mathbb{P}}\left[#1\right]}

\newcommand{\expn}[1]{\mathrm{\mathbb{E}}\left[#1\right]}

\def\COMMENT#1{}
\def\TASK#1{}
\let\TASK=\footnote             % COMMENT OUT for clean output
\begin{document}

\title{A hitting-time version of Kotzig's conjecture in random graphs}

\author{
Stefan Glock \thanks{Fakultät für Informatik und Mathematik, Universität Passau, Germany.
\emph{Email}: \href{mailto:stefan.glock@uni-passau.de}{\tt stefan.glock@uni-passau.de}, \href{mailto:amedeo.sgueglia@uni-passau.de}{\tt amedeo.sgueglia@uni-passau.de}.}
\and
Lyuben Lichev \thanks{Institute of Statistics and Mathematical Methods in Economics, TU Wien, A-1040 Vienna, Austria.
\emph{Email}: \href{mailto:lyuben.lichev@tuwien.ac.at}{\tt lyuben.lichev@tuwien.ac.at}}
\and
Amedeo Sgueglia \footnotemark[1]
}

\date{}

\maketitle

\begin{abstract} 
A family of perfect matchings on the same vertex set is said to have the Kotzig property if the union of any pair of matchings forms a Hamilton cycle.
In 1963, Anton Kotzig famously conjectured that, when $n$ is even, the complete graph on $n$ vertices can be decomposed into $n-1$ perfect matchings with the Kotzig property.
A closely related conjecture of Wormald states that, for fixed $d \ge 3$, with high probability a random $d$-regular graph contains $d$ perfect matchings with the Kotzig property. 
Recently, the first and third authors considered a version of these problems for the binomial random graph and proved that, for every fixed $k$, there is $C=C(k)$ such that, with high probability, the binomial random graph $G(n,C\log n/n)$ contains $k$ perfect matchings with the Kotzig property.
We strengthen this result by showing a hitting time version:
for every fixed $k$, with high probability, the hitting time in the random graph process for finding $k$ edge-disjoint perfect matchings with the Kotzig property is the same as the hitting time for a graph with minimum degree $k$.
This is a generalisation of the celebrated result of Bollob\'as and Frieze for even $n$ establishing the hitting time for $k$ edge-disjoint Hamilton cycles.

The proof uses several ideas.
We first design a coupling, via a randomised algorithm, of the random graph at the hitting time with $k$ almost-independent perfect matchings and a set of random reservoir edges, which may be of independent interest. 
We then prove various statistical properties of these matchings, in particular concerning the cycle structure of the $2$-factors formed by pairs of matchings. This includes a novel concept of ``shuffledness'', which we establish through a careful exposure argument and a concentration inequality for random permutations. Finally, we show that the initial $2$-factors can be merged into Hamilton cycles, using the random reservoir edges, by analysing certain P\'osa-type rotations via branching processes.
\end{abstract}

\section{Introduction}
\label{sec:intro}

A \defn{perfect matching} of a graph is a collection of pairwise vertex-disjoint edges covering all vertices.
For a graph to have a perfect matching, the number of its vertices must be even and thus we always assume this is the case. 
A partition of the edge set of a graph into perfect matchings is called a \emph{$1$-factorisation}.
The following standard construction (depicted for $n=6$ but readily generalisable for all even integers) shows that, when $n$ is even, the complete graph $K_n$ on $n$ vertices has a $1$-factorisation. 

\begin{figure}[htbp]
  \centering
  \begin{tikzpicture}[scale=0.8, every node/.style={draw, circle, fill=black, inner sep=1.5pt}]
    
    % ==================== CONFIGURATION F_1 ====================
    \begin{scope}[xshift=0.0cm]
      \coordinate (O) at (0,0);
      \coordinate (V1) at (90:1.2);
      \coordinate (V2) at (162:1.2);
      \coordinate (V3) at (234:1.2);
      \coordinate (V4) at (306:1.2);
      \coordinate (V5) at (378:1.2);
      
      % Background edges of K_6
      \draw[lightgray, thin] (O) -- (V1);
      \draw[lightgray, thin] (V1) -- (V2);
      \draw[lightgray, thin] (V1) -- (V3);
      \draw[lightgray, thin] (V1) -- (V4);
      \draw[lightgray, thin] (V1) -- (V5);
      \draw[lightgray, thin] (O) -- (V2);
      \draw[lightgray, thin] (V2) -- (V3);
      \draw[lightgray, thin] (V2) -- (V4);
      \draw[lightgray, thin] (V2) -- (V5);
      \draw[lightgray, thin] (O) -- (V3);
      \draw[lightgray, thin] (V3) -- (V4);
      \draw[lightgray, thin] (V3) -- (V5);
      \draw[lightgray, thin] (O) -- (V4);
      \draw[lightgray, thin] (V4) -- (V5);
      \draw[lightgray, thin] (O) -- (V5);
      
      % 1-factor edges for F_1
      \draw[ultra thick, red] (O) -- (V1);
      \draw[ultra thick, red] (V2) -- (V5);
      \draw[ultra thick, red] (V3) -- (V4);
      
      % Vertices
      \node at (O) {};
      \node at (V1) {};
      \node at (V2) {};
      \node at (V3) {};
      \node at (V4) {};
      \node at (V5) {};
    \end{scope}
    
    % ==================== CONFIGURATION F_2 ====================
    \begin{scope}[xshift=3.6cm]
      \coordinate (O) at (0,0);
      \coordinate (V1) at (90:1.2);
      \coordinate (V2) at (162:1.2);
      \coordinate (V3) at (234:1.2);
      \coordinate (V4) at (306:1.2);
      \coordinate (V5) at (378:1.2);
      
      % Background edges of K_6
      \draw[lightgray, thin] (O) -- (V1);
      \draw[lightgray, thin] (V1) -- (V2);
      \draw[lightgray, thin] (V1) -- (V3);
      \draw[lightgray, thin] (V1) -- (V4);
      \draw[lightgray, thin] (V1) -- (V5);
      \draw[lightgray, thin] (O) -- (V2);
      \draw[lightgray, thin] (V2) -- (V3);
      \draw[lightgray, thin] (V2) -- (V4);
      \draw[lightgray, thin] (V2) -- (V5);
      \draw[lightgray, thin] (O) -- (V3);
      \draw[lightgray, thin] (V3) -- (V4);
      \draw[lightgray, thin] (V3) -- (V5);
      \draw[lightgray, thin] (O) -- (V4);
      \draw[lightgray, thin] (V4) -- (V5);
      \draw[lightgray, thin] (O) -- (V5);
      
      % 1-factor edges for F_2
      \draw[ultra thick, red] (O) -- (V2);
      \draw[ultra thick, red] (V3) -- (V1);
      \draw[ultra thick, red] (V4) -- (V5);
      
      % Vertices
      \node at (O) {};
      \node at (V1) {};
      \node at (V2) {};
      \node at (V3) {};
      \node at (V4) {};
      \node at (V5) {};
    \end{scope}
    
    % ==================== CONFIGURATION F_3 ====================
    \begin{scope}[xshift=7.2cm]
      \coordinate (O) at (0,0);
      \coordinate (V1) at (90:1.2);
      \coordinate (V2) at (162:1.2);
      \coordinate (V3) at (234:1.2);
      \coordinate (V4) at (306:1.2);
      \coordinate (V5) at (378:1.2);
      
      % Background edges of K_6
      \draw[lightgray, thin] (O) -- (V1);
      \draw[lightgray, thin] (V1) -- (V2);
      \draw[lightgray, thin] (V1) -- (V3);
      \draw[lightgray, thin] (V1) -- (V4);
      \draw[lightgray, thin] (V1) -- (V5);
      \draw[lightgray, thin] (O) -- (V2);
      \draw[lightgray, thin] (V2) -- (V3);
      \draw[lightgray, thin] (V2) -- (V4);
      \draw[lightgray, thin] (V2) -- (V5);
      \draw[lightgray, thin] (O) -- (V3);
      \draw[lightgray, thin] (V3) -- (V4);
      \draw[lightgray, thin] (V3) -- (V5);
      \draw[lightgray, thin] (O) -- (V4);
      \draw[lightgray, thin] (V4) -- (V5);
      \draw[lightgray, thin] (O) -- (V5);
      
      % 1-factor edges for F_3
      \draw[ultra thick, red] (O) -- (V3);
      \draw[ultra thick, red] (V4) -- (V2);
      \draw[ultra thick, red] (V5) -- (V1);
      
      % Vertices
      \node at (O) {};
      \node at (V1) {};
      \node at (V2) {};
      \node at (V3) {};
      \node at (V4) {};
      \node at (V5) {};
    \end{scope}
    
    % ==================== CONFIGURATION F_4 ====================
    \begin{scope}[xshift=10.8cm]
      \coordinate (O) at (0,0);
      \coordinate (V1) at (90:1.2);
      \coordinate (V2) at (162:1.2);
      \coordinate (V3) at (234:1.2);
      \coordinate (V4) at (306:1.2);
      \coordinate (V5) at (378:1.2);
      
      % Background edges of K_6
      \draw[lightgray, thin] (O) -- (V1);
      \draw[lightgray, thin] (V1) -- (V2);
      \draw[lightgray, thin] (V1) -- (V3);
      \draw[lightgray, thin] (V1) -- (V4);
      \draw[lightgray, thin] (V1) -- (V5);
      \draw[lightgray, thin] (O) -- (V2);
      \draw[lightgray, thin] (V2) -- (V3);
      \draw[lightgray, thin] (V2) -- (V4);
      \draw[lightgray, thin] (V2) -- (V5);
      \draw[lightgray, thin] (O) -- (V3);
      \draw[lightgray, thin] (V3) -- (V4);
      \draw[lightgray, thin] (V3) -- (V5);
      \draw[lightgray, thin] (O) -- (V4);
      \draw[lightgray, thin] (V4) -- (V5);
      \draw[lightgray, thin] (O) -- (V5);
      
      % 1-factor edges for F_4
      \draw[ultra thick, red] (O) -- (V4);
      \draw[ultra thick, red] (V5) -- (V3);
      \draw[ultra thick, red] (V1) -- (V2);
      
      % Vertices
      \node at (O) {};
      \node at (V1) {};
      \node at (V2) {};
      \node at (V3) {};
      \node at (V4) {};
      \node at (V5) {};
    \end{scope}
    
    % ==================== CONFIGURATION F_5 ====================
    \begin{scope}[xshift=14.4cm]
      \coordinate (O) at (0,0);
      \coordinate (V1) at (90:1.2);
      \coordinate (V2) at (162:1.2);
      \coordinate (V3) at (234:1.2);
      \coordinate (V4) at (306:1.2);
      \coordinate (V5) at (378:1.2);
      
      % Background edges of K_6
      \draw[lightgray, thin] (O) -- (V1);
      \draw[lightgray, thin] (V1) -- (V2);
      \draw[lightgray, thin] (V1) -- (V3);
      \draw[lightgray, thin] (V1) -- (V4);
      \draw[lightgray, thin] (V1) -- (V5);
      \draw[lightgray, thin] (O) -- (V2);
      \draw[lightgray, thin] (V2) -- (V3);
      \draw[lightgray, thin] (V2) -- (V4);
      \draw[lightgray, thin] (V2) -- (V5);
      \draw[lightgray, thin] (O) -- (V3);
      \draw[lightgray, thin] (V3) -- (V4);
      \draw[lightgray, thin] (V3) -- (V5);
      \draw[lightgray, thin] (O) -- (V4);
      \draw[lightgray, thin] (V4) -- (V5);
      \draw[lightgray, thin] (O) -- (V5);
      
      % 1-factor edges for F_5
      \draw[ultra thick, red] (O) -- (V5);
      \draw[ultra thick, red] (V1) -- (V4);
      \draw[ultra thick, red] (V2) -- (V3);
      
      % Vertices
      \node at (O) {};
      \node at (V1) {};
      \node at (V2) {};
      \node at (V3) {};
      \node at (V4) {};
      \node at (V5) {};
    \end{scope}
    
  \end{tikzpicture}
  \caption{A $1$-factorisation of $K_6$.}
\end{figure}
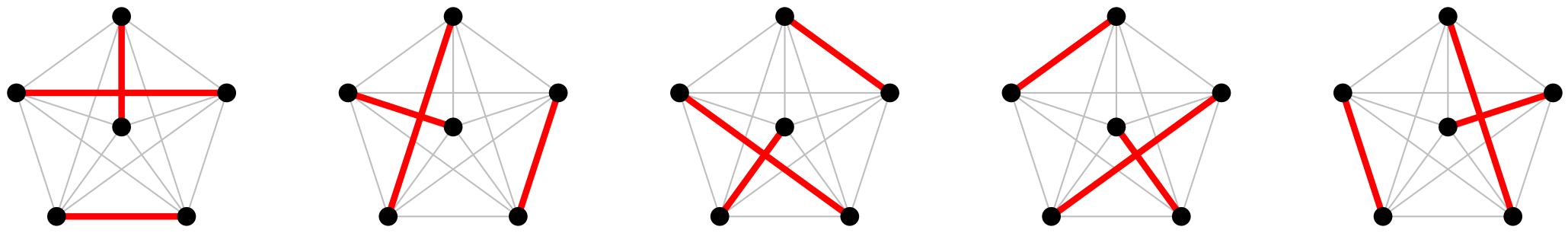

Kotzig observed that, if $n-1$ is prime, then each pair of perfect matchings in the construction above induces a Hamilton cycle: we call this property of a family of perfect matchings the \emph{Kotzig property}.
A $1$-factorisation having the Kotzig property is called \defn{perfect}.
Based on his observation, Kotzig made the following conjecture at the first international conference devoted to graph theory in 1963 (see~\cite{KO:14c}), nowadays known as \emph{Kotzig's conjecture} or \defn{perfect $1$-factorisation conjecture}~\cite{kotzig:64}.

\begin{conj}\label{conj:kotzig}
For every even integer $n \ge 4$, $K_n$ has a perfect $1$-factorisation.
\end{conj}

Apart from Kotzig's observation, the conjecture was proved by Anderson~\cite{anderson:73} when $n/2$ is prime, and is known to hold for $n\le 62$ and a few more sporadic values~\cite{Pik19,Ros19}.
Moreover, Cheng and the third author~\cite{CS:26+} proved that it holds asymptotically: namely, for even $n$, the complete graph $K_n$ contains a $1$-factorisation where $n-o(n)$ matchings have the Kotzig property.
For a rather complete account on the progress on Kotzig's conjecture and connections to the acyclic edge-chromatic number, to Latin squares and to orthogonal 1-factorisations, see the discussion in~\cite{GS:25} and the survey paper~\cite{Ros19}.

A related conjecture, which is also still open, is due to Wormald~\cite{wormald:99}.

\begin{conj}\label{conj:wormald}
    For fixed $d \ge 3$, a random $d$-regular graph whp~has a perfect $1$-factorisation.
\end{conj}

Recently, the first and third authors~\cite{GS:25} considered a variant of these problems for binomial random graphs. 
They showed that, for every fixed $k$, there is $C=C(k)$ such that, with high probability\footnote{That is, with probability tending to 1 as $n\to\infty$.}, 
the binomial random graph $G(n,C\log n/n)$ contains a set of $k$ perfect matchings with the Kotzig property.
In fact, their main result answers the question of counting, given any $k$ perfect matchings of $K_n$ with small overlap, the number of perfect matchings of $K_n$ which form a Hamilton cycle with all the given ones. The result on the binomial random graph then follows via the breakthroughs on the expectation threshold conjecture.
While this is sharp up to the value of the constant $C$, as $\log n/n$ is a sharp threshold for connectivity in $G(n,p)$, the only obvious necessary condition for the binomial random graph to contain $k$ such perfect matchings is that its minimum degree is at least~$k$.
Glock and Sgueglia~\cite{GS:25} suggested that in the random graph process, 
with high probability, $k$ perfect matchings with the Kotzig property appear at the moment when the minimum degree of the graph becomes exactly~$k$. 
In fact, many similar hitting time results for `sparse' graph properties hold in the binomial random graph; examples include the existence of a perfect matching~\cite{BT85,ER66}, Hamiltonicity~\cite{AKS:85,bollobas:84}, the existence of $k$ edge-disjoint Hamilton cycles~\cite{BF:85}, $k$-connectivity of the $k$-core~\cite{Luc91}, and similar phenomena also hold in other random graph models~\cite{BBKMW11,Bol90,CE-DGKO24,IL25,MPW11,Pen03}.

In this work, we confirm the hitting time result suggested in~\cite{GS:25}.
To state our main result, we first define the random graph process on $n$ vertices: it starts with the empty graph $G(0)$ and, for every $1\le t\le \binom{n}{2}$, one of the edges missing in $G(t-1)$ is chosen uniformly at random and is added to $G(t-1)$ to form $G(t)$.
For a monotone increasing graph property, its \emph{hitting time} is the smallest~$t$ such that $G(t)$ has this property.
We denote by $\tau_k$ the hitting time for the property of having minimum degree at least $k$.

\begin{theorem}\label{thm:main}
For every fixed $k$, with high probability, $G(\tau_k)$ contains $k$ edge-disjoint perfect matchings with the Kotzig property.
\end{theorem}

Note that this is a generalisation of the celebrated result of Bollob\'as and Frieze~\cite{BF:85} (for even $n$) which established that with high probability $G(\tau_k)$ contains $\lfloor k/2\rfloor$ disjoint Hamilton cycles (plus a disjoint perfect matching if $k$ is odd).
Indeed, these Hamilton cycles can be obtained by arbitrarily pairing the matchings given by Theorem~\ref{thm:main} and taking their union.
Moreover, since it is well-known (see e.g.\ Theorem~4.3 and its proof in~\cite{FK:15}) that
\begin{equation}\label{eq:tauk}
\text{with high probability } \tau_k=\frac{(\log n + (k-1) \log \log n\pm \omega(1))n}{2},
\end{equation}
Theorem~\ref{thm:main} also implies that, for every fixed $k$ and $p=p(n)$ with $np-\log n - (k-1) \log \log n = \omega(1)$, with high probability, the binomial random graph $G(n,p)$ contains $k$ perfect matchings having the Kotzig property.
Finally, we remark that the hitting time results for perfect matching and Hamilton cycle from~\cite{AKS:85,bollobas:84,BT85,ER66} already show Theorem~\ref{thm:main} for $k\in \{1,2\}$. 

\paragraph{Plan of the paper.} The remainder of our work is organised as follows. We finish the current section with a proof outline, and include notation and several preliminary results in Section~\ref{sec:prelims}.
Sections~\ref{sec:coupling} and~\ref{sec:merging} are the heart of our proof of Theorem~\ref{thm:main}, and we finish with some concluding remarks and open problems in Section~\ref{sec:conclusion}.

\subsection{Proof outline}\label{subsec:outline}
We now give an overview of our proof strategy, which is very different from the one in~\cite{GS:25} and is based on P\'osa-style rotations.
We first explain a simplified setup which can be used to give an alternative proof of the random graph result in~\cite{GS:25}.
Then, we will highlight the difficulties in adapting the approach to the hitting time setting and outline how we address them.

Suppose we want to show that, for fixed $k$, if $C=C(k)$ is large enough and $p \ge C (\log n)/n$, then with high probability $G \sim G(n,p)$ contains $k$ perfect matchings with the Kotzig property.
Let $G_1, \dots, G_{k+1}$ be $k+1$ mutually independent random graphs distributed as $G(n,\frac{p}{k+1})$. 
Then, their union is distributed as a binomial random graph with edge probability 
\[1-\left(1-\frac{p}{k+1}\right)^{k+1} \le \frac{p}{k+1}(k+1)= p\] 
and thus we can assume $G(n,p) \supseteq G_1 \cup \dots \cup G_{k+1}$.
Since $p/(k+1)$ is above the threshold for the existence of a perfect matching, with high probability we can find $k$ independent and uniformly random perfect matchings $M_1, \dots, M_k$ with $M_i$ being a perfect matching of~$G_i$. Each of the $2$-factors $M_i\cup M_j$ will then only consist of few cycles.

Our aim is to use the edges of $G_{k+1}$ (which can be thought of as a \emph{random reservoir}) to consecutively merge cycles of the $2$-factors until they become Hamilton cycles.
Suppose we want to merge $M_i \cup M_j$ into a Hamilton cycle by replacing some of the edges of $M_i$ with some of $G_{k+1}$.
A naive approach would be to choose two edges $ab,a'b'$ of $M_i$ on different cycles and then replace them in $M_i$ with the two edges $ab',a'b$. Then $M_i$ would still be a perfect matching, and the number of cycles of $M_i\cup M_j$ would be reduced by one. 
This strategy clearly comes with a lot of caveats. First of all, this operation is only possible when $ab',a'b$ are also edges of $G_{k+1}$, and such a pair will in fact most likely not exist. Secondly, even if the two edges were present, it could be that performing this switch would 
increase the number of cycles formed by $M_i$ and another matching.

To overcome these issues, we follow an approach which is inspired by P\'osa's rotation method. 
In fact, the first insight of our work is 
to find a way to perform P\'osa rotations while controlling the cycle structure of the union of \emph{any} pair of matchings.
Let $e=uv$ be an edge of $M_i$, delete it and let $P$ be the path of $M_i \cup M_j \setminus \{e\}$ with endpoints
$u$ and $v$.
Now, consider all neighbours of $v$ in~$G_{k+1}$.
Say that $x$ is such a neighbour and let $v_1$ be the vertex that $x$ is matched to by~$M_i$. Suppose that we decide to replace the edge $xv_1$ with $xv$ in $M_i$.
Then three scenarios are possible (see~Figure~\ref{fig:3cases}):
\begin{enumerate}[label=(\alph*)]
    \item $x$ does not belong to $P$ and thus it belongs to another cycle $C$ of $M_i \cup M_j$, in which case $P$ gets merged with the path $C\setminus \{xv_1\}$;
    \item $x$ belongs to $P$ and, when traversing $P$ from $u$ to $v$, $x$ appears before $v_1$, in which case $P$ gets rotated and the number of cycles is preserved;
    \item $x$ belongs to $P$ and, when traversing $P$ from $u$ to $v$, $v_1$ appears before $x$, in which case $P$ splits into a new cycle and a path.
\end{enumerate}
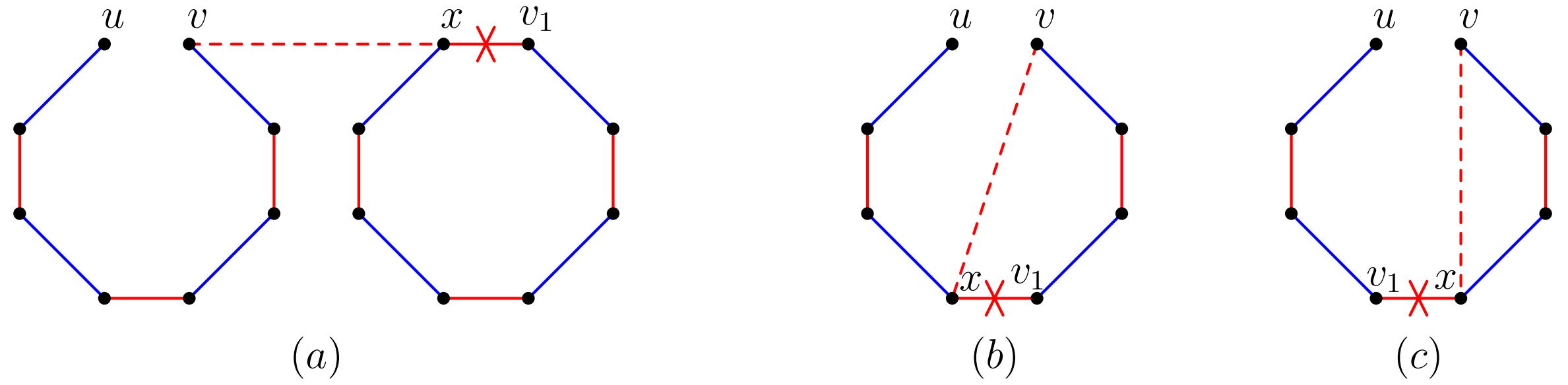
\begin{figure}[htp]
    \centering
    \begin{tikzpicture}[scale=0.85,line cap=round,line join=round,x=1cm,y=1cm]
\clip(-9.15,1) rectangle (9.15,5.5);
\draw [line width=0.8pt,color=blue] (-8,5)-- (-9,4);
\draw [line width=0.8pt,color=blue] (-9,3)-- (-8,2);
\draw [line width=0.8pt,color=blue] (-7,2)-- (-6,3);
\draw [line width=0.8pt,color=blue] (-7,5)-- (-6,4);
\draw [line width=0.8pt,color=blue] (-4,5)-- (-5,4);
\draw [line width=0.8pt,color=red] (-5,4)-- (-5,3);
\draw [line width=0.8pt,color=blue] (-5,3)-- (-4,2);
\draw [line width=0.8pt,color=red] (-4,2)-- (-3,2);
\draw [line width=0.8pt,color=blue] (-3,2)-- (-2,3);
\draw [line width=0.8pt,color=red] (-2,3)-- (-2,4);
\draw [line width=0.8pt,color=blue] (-2,4)-- (-3,5);
\draw [line width=0.8pt,color=red] (-3,5)-- (-4,5);
\draw [line width=0.8pt,color=red] (-9,4)-- (-9,3);
\draw [line width=0.8pt,color=red] (-8,2)-- (-7,2);
\draw [line width=0.8pt,color=red] (-6,3)-- (-6,4);
\draw [line width=0.8pt,dash pattern=on 3pt off 3pt,color=red] (-7,5)-- (-4,5);
\draw [line width=0.8pt,color=red] (-3.6,5.2)-- (-3.4,4.8);
\draw [line width=0.8pt,color=red] (-3.4,5.2)-- (-3.6,4.8);
\draw [line width=0.8pt,color=blue] (2,5)-- (1,4);
\draw [line width=0.8pt,color=red] (1,4)-- (1,3);
\draw [line width=0.8pt,color=blue] (1,3)-- (2,2);
\draw [line width=0.8pt,color=red] (2,2)-- (3,2);
\draw [line width=0.8pt,color=blue] (3,2)-- (4,3);
\draw [line width=0.8pt,color=red] (4,3)-- (4,4);
\draw [line width=0.8pt,color=blue] (4,4)-- (3,5);
\draw [line width=0.8pt,dash pattern=on 3pt off 3pt,color=red] (3,5)-- (2,2);
\draw [line width=0.8pt,color=blue] (7,5)-- (6,4);
\draw [line width=0.8pt,color=red] (6,4)-- (6,3);
\draw [line width=0.8pt,color=blue] (6,3)-- (7,2);
\draw [line width=0.8pt,color=red] (7,2)-- (8,2);
\draw [line width=0.8pt,color=blue] (8,2)-- (9,3);
\draw [line width=0.8pt,color=red] (9,3)-- (9,4);
\draw [line width=0.8pt,color=blue] (9,4)-- (8,5);
\draw [line width=0.8pt,dash pattern=on 3pt off 3pt,color=red] (8,5)-- (8,2);
\draw [line width=0.8pt,color=red] (2.4,2.2)-- (2.6,1.8);
\draw [line width=0.8pt,color=red] (2.6,2.2)-- (2.4,1.8);
\draw [line width=0.8pt,color=red] (7.4,2.2)-- (7.6,1.8);
\draw [line width=0.8pt,color=red] (7.6,2.2)-- (7.4,1.8);
\begin{scriptsize}
\draw [fill=black] (-8,5) circle (1.9pt);
\draw[color=black] (-7.9,5.3) node {\large{$u$}};
\draw [fill=black] (-9,4) circle (1.9pt);
\draw [fill=black] (-9,3) circle (1.9pt);
\draw [fill=black] (-8,2) circle (1.9pt);
\draw [fill=black] (-7,2) circle (1.9pt);
\draw [fill=black] (-6,3) circle (1.9pt);
\draw [fill=black] (-7,5) circle (1.9pt);
\draw[color=black] (-6.9,5.3) node {\large{$v$}};
\draw [fill=black] (-6,4) circle (1.9pt);
\draw [fill=black] (-4,5) circle (1.9pt);
\draw[color=black] (-3.9,5.3) node {\large{$x$}};
\draw [fill=black] (-5,4) circle (1.9pt);
\draw [fill=black] (-5,3) circle (1.9pt);
\draw [fill=black] (-4,2) circle (1.9pt);
\draw[color=black] (-5.5,1.3) node {\large{$(a)$}};
\draw [fill=black] (-3,2) circle (1.9pt);
\draw [fill=black] (-2,3) circle (1.9pt);
\draw [fill=black] (-2,4) circle (1.9pt);
\draw [fill=black] (-3,5) circle (1.9pt);
\draw[color=black] (-2.9,5.3) node {\large{$v_1$}};
\draw [fill=black] (2,5) circle (1.9pt);
\draw[color=black] (2.1,5.3) node {\large{$u$}};
\draw [fill=black] (1,4) circle (1.9pt);
\draw [fill=black] (1,3) circle (1.9pt);
\draw [fill=black] (2,2) circle (1.9pt);
\draw[color=black] (2.22,2.2) node {\large{$x$}};
\draw [fill=black] (3,2) circle (1.9pt);
\draw[color=black] (2.9,2.25) node {\large{$v_1$}};
\draw[color=black] (2.5,1.3) node {\large{$(b)$}};
\draw [fill=black] (4,3) circle (1.9pt);
\draw [fill=black] (4,4) circle (1.9pt);
\draw [fill=black] (3,5) circle (1.9pt);
\draw[color=black] (3.1,5.3) node {\large{$v$}};
\draw [fill=black] (7,5) circle (1.9pt);
\draw[color=black] (7.1,5.3) node {\large{$u$}};
\draw [fill=black] (6,4) circle (1.9pt);
\draw [fill=black] (6,3) circle (1.9pt);
\draw [fill=black] (7,2) circle (1.9pt);
\draw[color=black] (7.1,2.22) node {\large{$v_1$}};
\draw [fill=black] (8,2) circle (1.9pt);
\draw[color=black] (7.82,2.22) node {\large{$x$}};
\draw[color=black] (7.5,1.3) node {\large{$(c)$}};
\draw [fill=black] (9,3) circle (1.9pt);
\draw [fill=black] (9,4) circle (1.9pt);
\draw [fill=black] (8,5) circle (1.9pt);
\draw[color=black] (8.1,5.3) node {\large{$v$}};
\end{scriptsize}
\end{tikzpicture}
\caption{The edges of $M_i$ are depicted in red, those of $M_j$ in blue. The figure illustrates the possible effect of replacing the red edge $xv_1$ with a new edge $xv$.}
\label{fig:3cases}
\end{figure}
As we do not want to increase the number of cycles in any $2$-factor, we only allow the rotations which avoid the third scenario for each $j \in [k] \setminus \{i\}$ (see Definition~\ref{def:good}).

By performing one rotation or extension step as described above (assuming it is allowed), we obtain a new path with endpoints $u$ (as before) and $v_1$ (the new one).
Continuing in this way with~$v_1$ taking the role of $v$ and using the expansion properties of the random graph $G_{k+1}$, we may hope to find $n^{2/3}$ different ``alternative'' endpoints for~$v$ which can be reached using few ($\mathrm{polylog}(n)$) steps.
We then consider the original matching $M_i$ again and start a similar procedure from $u$, keeping $v$ as fixed endpoint, to find $n^{2/3}$ alternative endpoints for~$u$ (also reachable using few rotation/extension steps).
We show that, for at least $n^{2/3}/2$ alternative endpoints $v^*$ of the first process, there exist at least $n^{1/2}$ alternative endpoints $u^*$ of the second process which are compatible in the sense that we can apply the sequence of rotations/extensions from $v$ to $v^*$ and then that from $u$ to $u^*$ without creating any conflict. 
In particular, we get at least $n^{7/6}/2$ pairs $u^*v^*$ which are compatible and with high probability $G_{k+1}$ will contain one of these as an edge. We then modify $M_i$ as follows: we delete $uv$ from $M_i$, perform exactly those rotations/extensions that changed $v$ to $v^*$ and $u$ to $u^*$, and then add the edge $u^*v^*$.
If, at the very beginning of the first procedure, we picked only neighbours of $v$ outside $P$ (which we can guarantee with high probability by taking the edge $e$ in a shortest cycle of $M_i \cup M_j$) and thus performed a rotation of type (a) in Figure~\ref{fig:3cases}, then the number of cycles in $M_i \cup M_j$ has decreased by at least one. 
Repeating this merging procedure until $M_i\cup M_j$ becomes a Hamilton cycle, consecutively for all pairs $i,j$, leads to the desired family of $k$ matchings with the Kotzig property.

A delicate argument in our strategy is showing that typically, for every edge $uv$, there is a large number of allowed rotations/extensions and that this property is stable under small modifications of the original perfect matchings, since we will have to perform several rotations.
Referring to the notation introduced above, suppose that we delete the edge $uv \in M_i$ and want to apply a sequence of rotations.
For convenience, we refer to the edges of $M_i$ as the red edges.
We encode the three scenarios (a)--(c) from Figure~\ref{fig:3cases} by giving to each red edge an orientation for each $j \in [k] \setminus \{i\}$ defined as follows:
if the edge belongs to the path of $M_i \cup M_j \setminus \{e\}$, we orient the edge as we traverse the path from $u$ to $v$, and otherwise, we orient it arbitrarily.
We call this the $j$-orientation of the edge and remark it depends on the edge $uv$ and changes after each rotation.
In particular, if we want to perform the rotation obtained by replacing the red edge $xv_1 \in M_i$ with the edge $xv$ of $R$ we require that each time the edge $xv_1$ belongs to the path of $M_i \cup M_j \setminus \{e\}$, its $j$-orientation is from $x$ to $v_1$.
We say that $xv_1$ is a good edge if the $j$-orientations for $j \in [k] \setminus \{i\}$ are all the same. 
We remark that we can allow the scenario~(a) to happen and thus, if the edge does not belong to the path of $M_i \cup M_j \setminus \{e\}$, its $j$-orientation does not really matter for deciding if we want to perform the rotation.
We nevertheless ask \emph{all} the $j$-orientations to be the same for simplicity of the presentation.

\begin{figure}[htp]
    \centering
    \begin{tikzpicture}[line cap=round,line join=round,x=1cm,y=1cm]
\clip(-9.65,2.3) rectangle (5.1,9.5);
\draw [line width=0.8pt,color=blue] (-7,9)-- (-8,8.5);
\draw [line width=0.8pt,color=blue] (-9,7.5)-- (-9.5,6.5);
\draw [line width=0.8pt,color=blue] (-9.5,5.5)-- (-9,4.5);
\draw [line width=0.8pt,color=blue] (-8,3.5)-- (-7,3);
\draw [line width=0.8pt,color=blue] (-6,3)-- (-5,3.5);
\draw [line width=0.8pt,color=blue] (-4,4.5)-- (-3.5,5.5);
\draw [line width=0.8pt,color=blue] (-3.5,6.5)-- (-4,7.5);
\draw [line width=0.8pt,color=blue] (-6,9)-- (-5,8.5);
\draw [line width=0.8pt,color=blue] (-1.5,6.5)-- (-1,7.5);
\draw [line width=0.8pt,color=blue] (0,8.5)-- (1,9);
\draw [line width=0.8pt,color=blue] (-1.5,5.5)-- (-1,4.5);
\draw [line width=0.8pt,color=blue] (0,3.5)-- (1,3);
\draw [line width=0.8pt,color=blue] (2,3)-- (3,3.5);
\draw [line width=0.8pt,color=blue] (4,4.5)-- (4.5,5.5);
\draw [line width=0.8pt,color=blue] (4.5,6.5)-- (4,7.5);
\draw [line width=0.8pt,color=blue] (3,8.5)-- (2,9);
\draw [line width=0.8pt,color=red] (-8,8.5)-- (-9,7.5);
\draw [line width=0.8pt,color=red] (-9.5,6.5)-- (-9.5,5.5);
\draw [line width=0.8pt,color=red] (-9,4.5)-- (-8,3.5);
\draw [line width=0.8pt,color=blue] (-8.358682476719926,3.8586824767199253) -- (-8.514131752328009,3.844550724391918);
\draw [line width=0.8pt,color=blue] (-8.358682476719926,3.8586824767199253) -- (-8.344550724391917,4.014131752328007);
\draw [line width=0.8pt,color=gray] (-8.5,4) -- (-8.655449275608083,3.9858682476719927);
\draw [line width=0.8pt,color=gray] (-8.5,4) -- (-8.485868247671993,4.155449275608082);
\draw [line width=0.8pt,color=red] (-7,3)-- (-6,3);
\draw [line width=0.8pt,color=blue] (-6.300146841976342,3) -- (-6.400073420988171,2.880088105185806);
\draw [line width=0.8pt,color=blue] (-6.300146841976342,3) -- (-6.400073420988171,3.1199118948141944);
\draw [line width=0.8pt,color=gray] (-6.5,3) -- (-6.599926579011829,2.880088105185806);
\draw [line width=0.8pt,color=gray] (-6.5,3) -- (-6.599926579011829,3.1199118948141944);
\draw [line width=0.8pt,color=red] (-5,3.5)-- (-4,4.5);
\draw [line width=0.8pt,color=red] (-3.5,5.5)-- (-3.5,6.5);
\draw [line width=0.8pt,color=blue] (-3.5,6.199853158023657) -- (-3.380088105185805,6.099926579011829);
\draw [line width=0.8pt,color=blue] (-3.5,6.199853158023657) -- (-3.6199118948141944,6.099926579011829);
\draw [line width=0.8pt,color=gray] (-3.5,6) -- (-3.380088105185805,5.900073420988171);
\draw [line width=0.8pt,color=gray] (-3.5,6) -- (-3.6199118948141944,5.900073420988171);
\draw [line width=0.8pt,color=red] (-4,7.5)-- (-5,8.5);
\draw [line width=0.8pt,color=red] (0,8.5)-- (-1,7.5);
\draw [line width=0.8pt,color=red] (-1.5,6.5)-- (-1.5,5.5);
\draw [line width=0.8pt,color=red] (3,3.5)-- (4,4.5);
\draw [line width=0.8pt,color=red] (4.5,5.5)-- (4.5,6.5);
\draw [line width=0.8pt,color=blue] (4.5,6.199853158023657-0.15) -- (4.619911894814195,6.300073420988171-0.15);
\draw [line width=0.8pt,color=blue] (4.5,6.199853158023657-0.15) -- (4.380088105185806,6.300073420988171-0.15);
\draw [line width=0.8pt,color=gray] (4.5,6) -- (4.619911894814195,5.900073420988171);
\draw [line width=0.8pt,color=gray] (4.5,6) -- (4.380088105185806,5.900073420988171);
\draw [line width=0.8pt,color=red] (4,7.5)-- (3,8.5);
\draw [line width=0.8pt,color=gray] (-7,9)-- (-9.5,6.5);
\draw [line width=0.8pt,color=gray] (-9.5,5.5)-- (-3.5,5.5);
\draw [line width=0.8pt,color=gray] (-3.5,6.5)-- (-9,4.5);
\draw [line width=0.8pt,color=gray] (-6,3)-- (-6,9);
\draw [line width=0.8pt,dash pattern=on 3pt off 3pt,color=red] (-9,4.5)-- (-6,9);
\draw [shift={(-6.994505225673301,4.2609895486533995)},line width=0.8pt,color=gray]  plot[domain=3.789451244611938:4.7080314981558296,variable=\t]({1*1.2610015203630822*cos(\t r)+0*1.2610015203630822*sin(\t r)},{0*1.2610015203630822*cos(\t r)+1*1.2610015203630822*sin(\t r)});
\draw [line width=0.8pt,color=red] (2,3)-- (1,3);
\draw [line width=0.8pt,color=blue] (1.3001468419763405,3) -- (1.4000734209881693,3.1199118948141944);
\draw [line width=0.8pt,color=blue] (1.3001468419763405,3) -- (1.4000734209881693,2.880088105185806);
\draw [line width=0.8pt,color=gray] (1.5,3) -- (1.599926579011827,3.1199118948141944);
\draw [line width=0.8pt,color=gray] (1.5,3) -- (1.599926579011827,2.880088105185806);
\draw [line width=0.8pt,color=red] (-1,4.5)-- (2,9);
\draw [line width=0.8pt,color=gray] (1,9)-- (-1.5,6.5);
\draw [line width=0.8pt,color=gray] (-1.5,5.5)-- (4.5,5.5);
\draw [line width=0.8pt,color=gray] (4.5,6.5)-- (-1,4.5);
\draw [line width=0.8pt,color=gray] (2,3)-- (2,9);
\draw [shift={(1.0013402051838907,4.252680410367781)},line width=0.8pt,color=gray]  plot[domain=3.7861636319694907:4.711319110798276,variable=\t]({1*1.2526811272902325*cos(\t r)+0*1.2526811272902325*sin(\t r)},{0*1.2526811272902325*cos(\t r)+1*1.2526811272902325*sin(\t r)});
\begin{scriptsize}
\draw [fill=black] (-7,9) circle (1.9pt);
\draw[color=black] (-6.9,9.3) node {\large{$u$}};
\draw [fill=black] (-8,8.5) circle (1.9pt);
\draw [fill=black] (-9,7.5) circle (1.9pt);
\draw [fill=black] (-9.5,6.5) circle (1.9pt);
\draw [fill=black] (-9.5,5.5) circle (1.9pt);
\draw [fill=black] (-9,4.5) circle (1.9pt);
\draw[color=black] (-8.98,4.8) node {\large{$x$}};
\draw [fill=black] (-8,3.5) circle (1.9pt);
\draw[color=black] (-7.9,3.7) node {\large{$v_1$}};
\draw [fill=black] (-7,3) circle (1.9pt);
\draw [fill=black] (-6,3) circle (1.9pt);
\draw [fill=black] (-5,3.5) circle (1.9pt);
\draw [fill=black] (-4,4.5) circle (1.9pt);
\draw [fill=black] (-3.5,5.5) circle (1.9pt);
\draw [fill=black] (-3.5,6.5) circle (1.9pt);
\draw [fill=black] (-4,7.5) circle (1.9pt);
\draw [fill=black] (-6,9) circle (1.9pt);
\draw[color=black] (-5.9,9.3) node {\large{$v$}};
\draw [fill=black] (-5,8.5) circle (1.9pt);
\draw [fill=black] (-1.5,6.5) circle (1.9pt);
\draw [fill=black] (-1,7.5) circle (1.9pt);
\draw [fill=black] (0,8.5) circle (1.9pt);
\draw [fill=black] (1,9) circle (1.9pt);
\draw[color=black] (1.1,9.3) node {\large{$u$}};
\draw [fill=black] (-1.5,5.5) circle (1.9pt);
\draw [fill=black] (-1,4.5) circle (1.9pt);
\draw[color=black] (-0.98,4.8) node {\large{$x$}};
\draw [fill=black] (0,3.5) circle (1.9pt);
\draw[color=black] (0.1,3.7) node {\large{$v_1$}};
\draw [fill=black] (1,3) circle (1.9pt);
\draw [fill=black] (2,3) circle (1.9pt);
\draw [fill=black] (3,3.5) circle (1.9pt);
\draw [fill=black] (4,4.5) circle (1.9pt);
\draw [fill=black] (4.5,5.5) circle (1.9pt);
\draw [fill=black] (4.5,6.5) circle (1.9pt);
\draw [fill=black] (4,7.5) circle (1.9pt);
\draw [fill=black] (3,8.5) circle (1.9pt);
\draw [fill=black] (2,9) circle (1.9pt);
\draw[color=black] (2.1,9.3) node {\large{$v$}};
\draw[color=black] (-6.5,2.6) node {\large{$f_2$}};
\draw[color=black] (-3.1,6) node {\large{$f_1$}};
\draw[color=black] (4.9,6) node {\large{$f_1$}};
\draw[color=black] (1.5,2.6) node {\large{$f_2$}};
\end{scriptsize}
\end{tikzpicture}
    \caption{Fix $k=3$. The figure on the left side shows some edges of three matchings depicted in three different colours, with red used for $M_i$, blue for $M_j$ and grey for $M_{j'}$ where $\{i,j,j'\} = [3]$.
    We delete the red edge $e=uv$ and observe that both $M_i \cup M_j \setminus \{e\}$ and $M_i \cup M_{j'} \setminus \{e\}$ have an alternating path with endpoints $u$ and $v$. 
    Moreover, the edge $xv_1$ is oriented from $x$ to $v_1$ for each of the colours $j,j'$ and thus we can perform the rotation replacing $xv_1$ with $xv$, whose outcome is depicted in the figure on the right side.
    Note the effect of the rotation on the orientations of the red edges $f_1$ and $f_2$.
    Before the rotation, all the orientations of $f_1$ were the same and the same was true for $f_2$.
    After the rotation, this is no longer the case for~$f_1$~while~it continues to be the case for $f_2$ where \emph{all} orientations flipped, maintaining the property.
    }
    \label{fig:flip}
\end{figure}
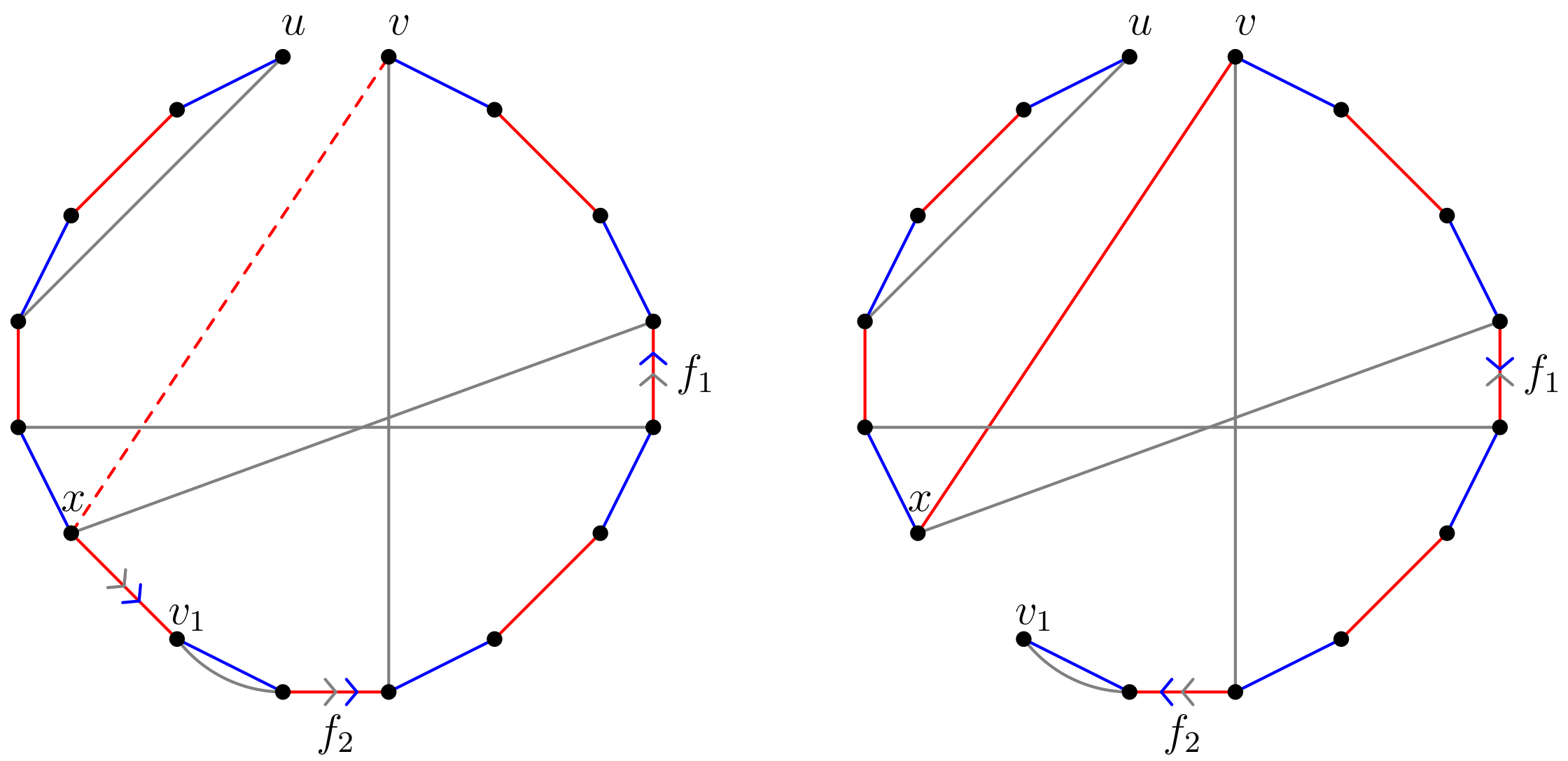

Heuristically, we expect that the probability that a red edge is good (for a fixed $e$) is $2 \cdot 2^{1-k} = 2^{2-k}$.
Even assuming one has concentration for the number of good edges, this alone is not sufficient as,
after performing a rotation, the orientations of the red edges change and an edge whose orientations were all the same may lose this property (see~Figure~\ref{fig:flip} for an example).
Instead we prove concentration for the number of good edges on segments:
every time we choose, for each $j \in [k] \setminus \{i\}$, any orientation of any sufficiently long alternating path $P_j$ of $M_i \cup M_j$, the (edge)-intersection $\cap_{j \in [k] \setminus \{i\}} P_j$ contains roughly the expected number of red edges with the same orientations.
This is captured by a novel concept of ``shuffledness'' (see~Definition~\ref{def:shuffled}).
The key property of shuffledness is its stability under small modifications of the initial matchings, which ensures the presence of many good edges at all times (with respect to the red edge we want to rotate; see~Proposition~\ref{prop:stays_shuffled}). We manage to reapply this argument until the matchings satisfy the Kotzig property.

The hitting time version of the result presents extra challenges and is even more involved.
In fact, we can still find independent subgraphs $G_1 \cup \dots \cup G_{k+1} \subseteq G(\tau_k)$ with $G_i \sim G(n,c(\log n)/n)$ for a constant $c>0$, but this time the constant $c$ cannot be chosen to be large enough to guarantee a perfect matching in the $G_i$.
A second insight of our work is an intricate coupling procedure (see Definition~\ref{def:multicore} and Theorem~\ref{thm:coupling_main}) that may be of independent interest. 
We show that with high probability the union of $G_1, \dots, G_k$ and a few more edges from $G(\tau_k)$ contains $k$ perfect matchings $M_1,\dots,M_k$ which are almost independent and almost uniformly distributed.
More precisely, there is a small set of vertices $B$ (informally, the vertices of small degree) such that the submatchings $M_1',\ldots,M_k'$ incident to $B$ are fixed somewhat greedily while $M_1\setminus M_1', \ldots, M_k\setminus M_k'$ are independent and uniformly random on their respective vertex sets.
Moreover, the procedure only reveals the edges of $G_1,\dots,G_k$ and the neighbourhood in $G(\tau_k)$ of the vertices in $B$; in particular, it does not reveal the edges of $G_{k+1}[V \setminus B]$.
Our construction ensures several properties of the (pairs of) matchings: namely, typically every 2-factor $M_i\cup M_j$ contains few cycles, the short cycles entirely avoid $B$ while long cycles contain long $B$-free segments, and the family is shuffled.
The graph $G_{k+1}[V \setminus B]$ is used to perform the rotations as above.
These properties allow the approach outlined in the simplified setup to be implemented, although being technically more demanding.

\section{Preliminaries}\label{sec:prelims}

We begin this section with a standard account of notation.

\vspace{1em}
\noindent
\textbf{General notation.} For a positive integer $n$, we write $[n] = \{1,\ldots,n\}$. Also, for real numbers $a,b,c$ with $c>0$, we write $a=b\pm c$ for the fact that $a\in [b-c,b+c]$.  
For two sets $A,B$, we denote by $A\Delta B$ their symmetric difference. We write $S_n$ for the set of permutations of $[n]$.
Throughout the paper, we systematically ignore rounding signs whenever this does not affect the validity of our arguments.\\
\textbf{Asymptotic notation.} We use standard asymptotic notation. When the implicit constants therein depend on some parameter, we indicate this parameter as a lower right index: for example, $O_c$ or~$\Omega_{\varepsilon}$.\\
\textbf{Probability.} For a random variable $X$ and a distribution $\cL$, we write $X\sim \cL$ to denote that $X$ is distributed according to $\cL$. 
The term \emph{with high probability} is often abbreviated \emph{whp} for convenience.\\
\textbf{Graphs.} For a graph $G$, we denote by $V(G)$ its vertex set and by $E(G)$ its edge set. We call $|V(G)|$ the \emph{order} of $G$ and $|E(G)|$ the \emph{size} of $G$.
For two sets of vertices $X,Y$ in a graph $G$, we write $E(X,Y)$ for the set of edges between $X$ and $Y$.
For a vertex $v \in V(G)$, we write $\deg(v,Y):=|E(\{v\},Y)|$.
We say that a matching $M$ \emph{saturates} a vertex $v$ if $v\in V(M)$.

\paragraph{Concentration inequalities.}

We first state the classic Chernoff bound, see e.g. \cite[Theorem~2.1]{JLR:00}.

\begin{lemma}\label{lem:chernoff}
For every binomial 
random variable $X$ and $t\ge 0$,
\begin{align*}
\prob{|X - \mathbb E[X]| \ge t} \le 2 \exp \left( - \frac {t^2}{2 (\mathbb E[X] + t/3)} \right).
\end{align*}
\end{lemma}

We next state the equally classic Azuma--Hoeffding concentration inequality for martingales, see e.g.\ \cite[Theorem~2.25]{JLR:00}.

\begin{lemma}\label{lem:azuma}
Consider a martingale $(X_i)_{i=0}^n$ with respect to a filtration $(\cF_i)_{i=0}^n$ where $|X_i - X_{i-1}| \le c_i$ almost surely for each $i\in [n]$. Then, for any $t > 0$,
\[
    \prob{|X_n - X_0| \ge t} \le 2 \exp\left( -\frac{t^2}{2 \sum_{i=1}^n c_i^2} \right).
\]
\end{lemma}

Finally, we show a version of the famous concentration inequality of McDiarmid~\cite{mcdiarmid:89}, also reminiscent of the concentration inequality for functions with bounded differences under switchings in random graphs~\cite[Theorem~2.19]{wormald:99}.
Given two permutations $\sigma,\sigma' \in S_n$, we say that they \defn{differ by a transposition} if there exists a transposition $\pi \in S_n$ such that $\sigma'=\sigma \circ \pi$.

\begin{lemma}
\label{lem:mcdiarmid_permutation}
Denote by $\cS:=S_n \times S_n \times \dots \times S_n $ the product of $k$ copies of $S_n$ and consider a function $f:\cS \to \bR$ such that, for all $i \in [k]$ and all permutations $\sigma_1,\dots,\sigma_k,\sigma_i' \in S_n$ such that $\sigma_i$ and $\sigma_i'$ differ by a transposition, we have that 
\[|f(\sigma_1,\ldots,\sigma_{i-1},\sigma_i,\sigma_{i+1},\ldots,\sigma_k)- f(\sigma_1,\ldots,\sigma_{i-1},\sigma_i',\sigma_{i+1},\ldots,\sigma_k)| \le c.\]
Then, if $S \in \cS$ is a uniformly random element of $\cS$, for every $t > 0$,
    \[ 
        \prob{|f(S)-\expn{f(S)}| \ge t} \le 2\exp\left(-\frac{t^2}{2nkc^2}\right)\, .
    \]
\end{lemma}

\begin{proof}
Fix $m=kn$.
We identify $S$ with a (random) string $X=(X_i)_{i=1}^m$ of length $m$ where every integer in $[n]$ is met exactly once in each of the intervals $(X_i)_{i=1}^n$, $(X_i)_{i=n+1}^{2n}$, \ldots, $(X_i)_{i=m-n+1}^m$. 
Define the filtration $(\mathcal{F}_i)_{i=0}^{m}$ generated by revealing $X$ one symbol at a time, and define the Doob martingale $M_i = \expn{f(S) \mid \mathcal{F}_i}$.
Note that $M_0=\expn{f(S)}$ and $M_m=f(S)$.

Upon revealing $X_1,\ldots,X_{i-1}$, denote by $R_i$ the set of available values for $X_i$ and, for every $x\in R_i$, define the family $T_x$ consisting of all realisations of the string $X$ given $\cF_{i-1}$ where $X_i=x$. 
Note that, for every $x,y \in R_i$, the transposition $(x \ y)$ between the $\lceil i/n\rceil$-th occurrences of $x$ and $y$ induces a bijection $\varphi_{x,y}$ between $T_x$ and $T_y$.
Thus, the triangle inequality and our hypothesis imply that 
\begin{align*}
\bigg|\expn{f(X)\mid \mathcal{F}_{i-1}, X_i = x} - \expn{f(X)\mid \mathcal{F}_{i-1}, X_i = y}\bigg| 
& =\bigg| \frac{1}{|T_x|} \sum_{Z \in T_x} f(Z) - \frac{1}{|T_y|} \sum_{Z' \in T_y} f(Z')\bigg| \\
&= \bigg| \frac{1}{|T_x|} \sum_{Z \in T_x} \Big(f(Z) - f(\varphi_{x,y}(Z)) \Big)\bigg| \\
&\le \frac{1}{|T_x|} \sum_{Z \in T_x} |f(Z) - f(\varphi_{x,y}(Z))| \le c\, .
\end{align*}
As a result, for every $i\in [m]$ and $x \in R_i$,
\begin{align*}
\bigg|\expn{f(X)\mid \mathcal{F}_{i-1}, X_i = x} - M_{i-1}\bigg|&= \bigg|\expn{f(X)\mid \mathcal{F}_{i-1}, X_i = x} - \frac{1}{|R_i|} \sum_{y\in R_i} \expn{f(X)\mid \mathcal{F}_{i-1}, X_i = y}\bigg|\\
&\le  \frac{1}{|R_i|} \cdot \sum_{y \in R_i} \bigg|\expn{f(X)\mid \mathcal{F}_{i-1}, X_i = x} - \expn{f(X)\mid \mathcal{F}_{i-1}, X_i = y}\bigg| \\
&\le c,
\end{align*}
where the first inequality follows from the triangle inequality.
Since the above is true for every $x \in R_i$, it follows that $|M_i-M_{i-1}| \le c$.
Since this is true for each $i \in [m]$, an application of Lemma~\ref{lem:azuma} finishes the proof.
\end{proof}

\paragraph{Coupling of $G(\tau_k)$.}
Recall that $\tau_k$ denotes the hitting time in the random graph process for the property of having minimum degree at least $k$.
We show the following simple coupling proposition.

\begin{prop}\label{prop:division}
Fix $k\ge 2$ and let $p_1=\frac{1}{k+1/4} \cdot \frac{\log n}{n}$ and $p_2=\frac{1}{8k+2} \cdot \frac{\log n}{n}$.
Consider $k+1$ independent random graphs $G_1,\ldots,G_{k+1}$ on a common vertex set $V$ with distribution $G_i \sim G(n,p_1)$ for $i \in [k]$ and $G_{k+1} \sim G(n,p_2)$. 
Then, there exists a coupling of these graphs and $G(\tau_k)$ such that whp
\[G_1\cup\ldots\cup G_{k+1}\subseteq G(\tau_k).\]
\end{prop}
\proof
The graph $G_1 \cup \dots \cup G_{k+1}$ can be coupled with the graph $G(n,p)$ where $1-p=(1-p_1)^k(1-p_2)$ since they have the same distribution.
Moreover 
\[p=1-(1-p_1)^k(1-p_2) = 1-(1-p_1)^k+p_2(1-p_1)^k \le kp_1+p_2 = \tfrac{8k+1}{8k+2} \cdot \tfrac{\log n}{n},\] 
where the inequality above is an application of Bernoulli's inequality. 
Let $N$ be the number of edges of $G(n,p)$ and observe $\expn{N} \le \frac{8k+1}{8k+2} \cdot \frac{n \log n}{2}$. 
Then,
\begin{equation*}
\begin{split}
\prob{N > \tau_k} \le \prob{N > (\log n)n/2}+\prob{\tau_k \le (\log n)n/2} = o(1)\, ,
\end{split}
\end{equation*}
where we used \eqref{eq:tauk}, the fact that $k\ge 2$, and the Chernoff bound.
As $G(n,p) = G(N)$ under the natural coupling with the random graph process, the statement follows.
\endproof

\paragraph{Perfect matchings in random graphs.}
We finish this section by showing that, in the binomial random graph $G(n,p)$ with $np=\Omega(\log n)$, every induced subgraph with minimum degree of order $\Omega(p n)$ and an even number of vertices contains a perfect matching.
The proof uses standard techniques but we include it for completeness. 
We recall the Tutte criterion for existence of a perfect matching: a graph $G$ contains a perfect matching if and only if, for every set of vertices $S \subseteq V(G)$, the number of connected components in $G\setminus S$ whose order is odd is at most $|S|$.

\begin{prop}\label{prop:min_degree_perf_matching}
Fix any $\eta , \delta> 0$ and $H \sim G(n,p)$ where $p \ge \delta \log n /n$.
Then, whp every induced subgraph $H' \subseteq H$ of even order and with minimum degree at least $\eta pn$ has a perfect matching.
\end{prop}
\begin{proof}
We set $\alpha_0 := \frac{\eta \e^{-16/(\eta \delta)}}{4 \e}$ and $m_0 := \frac{4 \log (np)}{p}$, and start by the following claim.
\begin{cl}\label{cl:H}
Whp the following properties hold for $H$:
\begin{itemize} 
    \item no set $X\subseteq V(H)$ of size $x \in [\eta pn /4, \alpha_0 n]$ satisfies $|E(H[X])| \ge \eta pnx/8$, where we implicitly assume we are in the range $p <4\alpha_0/\eta$, as otherwise the interval for $x$ is empty;
    \item the largest independent set of $H$ has size at most $m_0$;
    \item there is always an edge of $H$ between any two disjoint sets of size at least $\alpha_0 n$;
    \item for each $d \in [n]$, $s \le \min\{m_0,2d\}$ and any two disjoint sets $D$ and $S$ of size $d$ and $s$ respectively we have $|E(D,S)| < 3\eta pn d/4$. 
\end{itemize}
\end{cl}
\lateproof{Claim~\ref{cl:H}}
We start from the first point.
By a union bound over all sizes~$x$, the probability that the event from the claim fails is at most
\begin{align}
\sum_{x = \eta p n/4}^{\alpha_0 n} \binom{n}{x} \binom{\binom{x}{2}}{x\eta p n/8} p^{x\eta p n/8}\le \sum_{x = \eta p n/4}^{\alpha_0 n} \bigg(\frac{\e n}{x}\bigg)^{x} \bigg(\frac{\e p x^2/2}{x\eta p n/8}\bigg)^{x\eta p n/8}.\label{eq:ub}
\end{align}
Using that $x \le \alpha_0 n$ and the definition of $\alpha_0$ we have
\[ \frac{\e p x^2/2}{x\eta p n/8} = \frac{4 \e x }{\eta n} \le \frac{4 e (\alpha_0 n)}{\eta n} = \e^{-16/(\eta \delta)},\]
and thus the summand in \eqref{eq:ub} is at most
\[
\bigg(\frac{\e n}{x} \e^{-2p n/\delta}\bigg)^{x} \le \bigg(\frac{4\e n}{\eta \delta \log n} \e^{-2 \log n}\bigg)^{x} = \bigg(\frac{4\e}{\eta \delta \log n} \cdot \frac{1}{n}\bigg)^{x} \le \bigg(\frac{1}{n}\bigg)^x,\]
where we used that $p \ge \delta \log n/n$ and $x \ge \eta pn/4 \ge \eta \delta \log n/4$ for the first inequality. Therefore the sum in \eqref{eq:ub} is at most $\sum_{x = 1}^{\infty} (1/n)^{x} = o(1)$.

For the second point, note that the probability that there exists an independent set of size $m_0$ is at most
\[\binom{n}{m_0} (1-p)^{\binom{m_0}{2}} \le \left[ \frac{\e n}{m_0} \cdot \e^{-p(m_0 - 1)/2} \right]^{m_0} = o(1).\]

The third point holds as the probability that there exist two sets of size at least $\alpha_0 n$ without an edge between them is at most
\[ 2^{2n} \cdot (1-p)^{(\alpha_0 n)^2} \le 2^{2n} \cdot \e^{-\delta \alpha_0^2 n \log n} = o(1). \]

Finally, let $d \in [n]$ and $s \le \min\{m_0,2d\}$.
Then the probability that there exist disjoint sets $D$ and $S$ of size $d$ and $s$ with $|E(D,S)| \ge 3\eta pn d/4$ is at most 
\[\binom{n}{s} \binom{n}{d} \binom{ds}{ 3\eta pn d/4} p^{3\eta pn d/4} \le n^{s+d} \left( \frac{4 \e s}{3\eta n} \right)^{3\eta p n d/4}\le n^{3d} \left( \frac{4 \e m_0}{3\eta n} \right)^{3\eta p n d/4}.\]
Using that $(4em_0)/(3\eta n) = o(1)$, the latter expression is of order $o(n^{-2})$, and a union bound over all choices of $s$ and $d$ shows the last point of the claim.

As all properties fail with probability $o(1)$, the statement of the claim holds.
\endproof

From now on, we assume that $H$ satisfies the properties of Claim~\ref{cl:H}.
Assume for contradiction that there is an induced subgraph $H'$ with even order and minimum degree at least $\eta pn$ which contains no perfect matching. 
By the Tutte criterion, there exists a set $S \subseteq V(H')$ such that $H'\setminus S$ contains $q>s\coloneqq |S|$ components $C_1, \dots, C_q$ with odd orders. Note that $q \equiv s \pmod 2$ since $|V(H')|$ is even, implying that $q \ge s + 2$.

Fix a set $I$ containing one vertex from each component in $H'\setminus S$; in particular, $I$ is an independent set of size $q$.
By the second property of Claim~\ref{cl:H}, $|I| \le m_0$ and thus $q\le m_0$.

Next, recall that every vertex in an odd component $C_i$ has at least $\eta pn$ neighbours in $H'$, which are all contained in $C_i \cup S$. In particular, $2|E(C_i)| + |E(C_i, S)| \ge \eta pn |V(C_i)|$.
We classify the components $C_1,\ldots,C_q$ into three classes based on their orders: small ($|V(C_i)| \le \eta np/4$), medium ($\eta np/4 < |V(C_i)| \le \alpha_0 n$) and large ($|V(C_i)| > \alpha_0 n$).
Note that, when $p\ge 4\alpha_0/\eta$, the first and the last regime overlap while the middle regime is non-existent; in such cases, our argument apply to only the two non-trivial regimes.
We further denote by $q_s$, $q_m$ and $q_{\ell}$ the number of small, medium and large components, respectively, and define $D$ to be the union of (the vertex sets of) the small and the medium components and set $d:=|D|$.
Then, if $C_i$ is a small component, we have $|E(C_i)| \le |V(C_i)|^2/2 \le |V(C_i)|\eta np/8$ and the same bound is true for the middle components by the first property of Claim~\ref{cl:H}.
Therefore,
\begin{equation} \label{eq:e(D,S)}
\begin{split}
|E(D, S)| = \sum_{C_i \subseteq D} |E(C_i, S)| &\ge \sum_{C_i \subseteq D} \eta pn |V(C_i)| - 2|E(C_i)|\\
&\ge \sum_{C_i \subseteq D} \eta p n |V(C_i)| - 2|V(C_i)| \frac{\eta p n}{8}\ge  \frac{3\eta p n}{4} d.
\end{split}
\end{equation}
We consider two cases. If $s<3\eta np/4$, then the number of edges between $S$ and $D$ is at most $ds$, which contradicts the latter display unless $D=\emptyset$.
Hence, every component among $C_1,\ldots,C_q$ must be large. 
However, Tutte's condition requires $q \ge s + 2 \ge 2$ and, by the third property of Claim~\ref{cl:H}, there is an edge between any two components, a contradiction.

Next, suppose that $s \ge 3\eta np/4$. 
Since the number of large components is bounded by $1/\alpha_0$, the number of small and medium components satisfies $q_s+q_m \ge s - q_{\ell} \ge s/2$. 
As each component contains at least one vertex, $d \ge q_s+q_m \ge s/2$.
Therefore, since $s<q\le m_0$, we have $s \le \min\{m_0,2d\}$ and the fourth property of Claim~\ref{cl:H} contradicts~\eqref{eq:e(D,S)}.
\end{proof}

\section{Constructing the perfect matchings}
\label{sec:coupling}
Recall that $\tau_k$ denotes the hitting time in the random graph process for the property of having minimum degree at least $k$.
The aim of this section is to prove the following theorem.
\begin{theorem}
\label{thm:coupling_main}
    Fix $k \ge 2$. There is $\eps>0$ such that the following holds.
    Denote by $V$ the vertex set of the random graph process.
    Then, whp $G(\tau_k)$ can be coupled with $k$ perfect matchings $M_1, \dots, M_k$, a set of vertices $B$ and a set of edges $R$ such that $M_1 \cup \dots \cup M_k \cup R \subseteq G(\tau_k)$ and
    \begin{enumerate}[label=\upshape{\textbf{S\arabic*}}]
        \item \label{item:S1} $|B| \le n^{1-\eps}$;
        \item \label{item:S2} for each $i \in [k]$, if $M_i'$ denotes the edges of $M_i$ incident to vertices of $B$ and $M_i'':=M_i \setminus M_i'$, then $M_i''$ is a uniformly random perfect matching on $V\setminus{V(M_i')}$. Moreover, the matchings $M_1'',\dots,M_k''$ are mutually independent;
        \item \label{item:S3} the edges of $R$ appear independently with probability $p_2=\frac{1}{8k+2} \cdot \frac{\log n}{n}$ between vertices in $V\setminus B$.
    \end{enumerate}
\end{theorem}

In order to prove Theorem~\ref{thm:coupling_main}, we fix $k \ge 2$ and constants 
\[ 
    \delta_1 = \frac{1}{k+1/4}\, , \quad \delta_2=\frac{1}{8k+2} \quad \text{ and } \quad \eta\ll \eps\ll 1/k
\]
for the purposes of the entire section, and set 
\[ p_i=\delta_i (\log n)/n \text{ for } i=1,2, \quad p_H=1-(1-p_1)^k, \quad q = 2(\log n)/n, \quad \text{ and } \quad L = \eta \log n\, .\] 
Let $G_1, \dots, G_{k+1}$ be independent random graphs on the same vertex set~$V$ with distribution $G_i \sim G(n,p_1)$ for $i \in [k]$ and $G_{k+1} \sim G(n,p_2)$, and set $H:=G_1 \cup \dots \cup G_k$, noticing that it is distributed as $G(n,p_H)$.
Equip~$V$ with a fixed circular ordering. 

Since working with $G(\tau_k)$ is harder than working with binomial random graphs, we sandwich it as follows.
By~\eqref{eq:tauk} and Proposition~\ref{prop:division}, under the natural coupling of the graphs in the random graph process, whp $H \cup G_{k+1} \subseteq G(\tau_k) \subseteq G(n, q)$. We assume this relation in the remainder of the proof.
Observe that, for example, this allows us to work with $G(n,q)$ every time we want to upper bound the probability that $G(\tau_k)$ satisfies a certain increasing monotone property.

To prove Theorem~\ref{thm:coupling_main}, we use a randomised algorithm on the edges of $G(\tau_k)$ to output, with high probability, $k$ perfect matchings and a set $B$ with the desired properties.
The algorithm does not reveal the edges of $G_{k+1}$ in $V \setminus B$ and thus we will be able to take $R:=G_{k+1}[V \setminus B]$.
Its proof is a direct consequence of the lemmas developed in the next subsection and it is therefore omitted.

\subsection{Multi-core matching construction} 

\begin{defin}[Multi-core matching construction]\label{def:multicore}
The \emph{multi-core matching construction} is defined in two stages as follows.

\vspace{0.5em}
\noindent
\textbf{Stage 1.} 
Initiate a dynamic set $B\leftarrow \emptyset$ and dynamic matchings $M_1\leftarrow \emptyset,\ldots,M_k\leftarrow \emptyset$ which will grow over time.
Initiate also dynamic sets $V_1\gets V \setminus V(M_1)$, \dots, $V_k\gets V \setminus V(M_k)$.
We consecutively process the vertices in $V$ along the fixed cyclic order; note that one vertex may be processed several times.
To process a vertex $v\in V$ at a given step $i\in \mathbb N$, we check if the event 
\begin{equation}\label{eq:Bv}
\cB_v(i) := \{\text{there is } j\in [k] \text{ such that $\deg_{G_j}(v,V_j)<L$}\}
\end{equation} 
holds: we call such vertices \emph{bad} and gradually accumulate them in the set $B$. 
Note two crucial facts:
\begin{itemize}
\item since the matchings $M_1,\ldots,M_k$ increase throughout the process, for a fixed vertex $v$ which is processed at each of the steps $i_1 < i_2$, we have $\cB_{v}(i_1)\subseteq \cB_{v}(i_2)$, and
\item to confirm the event $\cB_v(i)$, we need to expose the indicator random variables of the events $\deg_{G_j}(v,V_j)<L$ for all $j\in [k]$ but not the precise degrees themselves.
\end{itemize}

Now, fix $v$ and $i$ as above. 
If $\cB_{v}(i)$ does not hold or $v\in B$, directly transition to the next vertex in $V$.
If $\cB_{v}(i)$ does hold and $v\notin B$, reveal all edges incident to $v$ in $G(\tau_k)[V\setminus B]$ and extend the matchings among $(M_j)_{j=1}^k$ which do not yet saturate $v$ by adding distinct edges saturating $v$ to each of them according to an arbitrary greedy deterministic rule. (Observe that the edge added to some $M_j$ belongs to $G(\tau_k)$ but not necessarily to $G_j$.)
If this task is impossible, we abort the construction and say that it \emph{terminates unsuccessfully}. 
Otherwise, update $B\leftarrow B\cup \{v\}$, $V_j\leftarrow V_j\setminus V(M_j)$ for every $j\in [k]$, and proceed to the next vertex.
Finally, transition to the second stage once $B$ remains unchanged for $n$ steps.

\vspace{0.5em}
\noindent
\textbf{Stage 2.} Assuming that the first stage did not terminate unsuccessfully, we start the second one with subgraphs $G_1'' := G_1[V_1],\ldots,G_k'' := G_k[V_k]$,
and observe that every vertex in $V\setminus B$ has at least $L$ neighbours in each of the graphs $G_1'',\ldots,G_k''$ it belongs to.
Sample independently uniformly random permutations $(\sigma_j)_{j=1}^k$ of $(V_j)_{j=1}^k$.
Then, for every $j\in [k]$, expose the graph $\sigma_j(G_j'')$ and find a perfect matching $M_j'''$ in it; if such is non-existent, we abort the construction and say that it \emph{terminates unsuccessfully}. 
Finally, for all $j\in [k]$, set $M_j''=\sigma_j^{-1}(M_j''')$, update $M_j\leftarrow M_j''\cup M_j$ and return $B,M_1,\ldots,M_k$. In this case, we say that the construction \emph{terminated successfully}.
\end{defin}

Until the end of this section, we show several key properties of the above construction.
First, we make an observation on the distribution of the graphs $G_1'',\dots,G_k''$.

\begin{lemma}\label{lem:unif}
Upon successful termination of Stage $1$, conditionally on the information revealed at Stage $1$ of the matching construction, the graphs $(G_i'')_{i=1}^k$ are jointly independent and their distributions are invariant under permutation of their respective vertex sets. 
In particular, upon successful termination of the entire matching construction, the matchings $(M_j'')_{j=1}^k$ are independent and uniformly random on their respective vertex sets.
\end{lemma}
\begin{proof}
Denote by $\tau_{\rm mc}$ the (random) number of steps in Stage 1 of the matching construction. 
Furthermore, for every $i\in [\tau_{\rm mc}]$, denote by $\cF_i$ the $\sigma$-algebra generated by the history of the process until step $i$ which includes:
\begin{itemize}
\item the indicators of the events $\deg_{G_j}(v,V_j)<L$ for all $j\in [k]$ and over all steps $i'\le i$ (where $v$ is the vertex processed at step $i'$), and
\item the indicator variables of the edges incident in $G(\tau_k)$ to the vertices entering $B$ over the first $i$ steps.
\end{itemize}
Moreover, for every $j\in [k]$ and for every vertex $v\in V$ processed at least once up to step $i$, denote by $V_j(v,i)$ the value of the set $V_j$ immediately after the last step in the interval $[i]$ when vertex $v$ was processed.

Recall that, for a vertex $v$ processed at each of the steps $i_1 < i_2$, we have $\cB_{v}(i_1)\subseteq \cB_{v}(i_2)$.
Using this fact and the definition of $\cF_i$, we show by induction over the steps of the matching construction that, for every step $i\in [\tau_{\rm mc}]$, conditionally on $\cF_i$, the random graphs $G_1,\ldots,G_k$ are jointly independent and their (conditional) distributions can be expressed as follows:
the distribution of the graph $G_j$ is the one of a binomial random graph with parameters $n$ and $p_1$ conditionally on: 
\begin{itemize}
\item the (indicators of the) edges incident in $G(\tau_k)$ to the set $B$, and
\item the intersection of the events $\deg_{G_j}(v,V_j(v,i))\ge L$ for all $v\in V\setminus B$ already processed at least once.
\end{itemize}
Note that, if $i' <i$, then $V_j(v,i') \supseteq V_j(v,i)$ but the inclusion can be strict as some more vertices may have left the set $V_j$ between steps $i'$ and $i$.

Since the graphs $G_1,\ldots,G_k$ are independent binomial random graphs, the base case is trivial.
Suppose that the induction hypothesis holds at all steps until $i-1$. 
If the currently processed vertex $v$ is already bad, its processing reveals no new information (so $\cF_i=\cF_{i-1}$) and the induction hypothesis is satisfied.
Otherwise, we expose the indicators of the events $\deg_{G_j}(v,V_j)<L$ for all $j\in [k]$ and consider two cases.

\vspace{0.5em}

\noindent
\textbf{Case 1.} Suppose that the event $\cB_i(v)$ fails.
Then, we update $V_j(v,i)$ and the only new information revealed is that $\deg_{G_j}(v,V_j(v,i))\ge L$ for all $j\in [k]$; note that the $j$-th of these events concerns only the random graph $G_j$ and none of the others.
Since the random graphs $G_1,\ldots,G_k$ were independent conditionally on $\cF_{i-1}$ at step $i-1$, the events $\deg_{G_j}(v,V_j(v,i))\ge L$ also are and hence the graphs $(G_j)_{j=1}^k$ remain independent conditionally on $\cF_i$ at step $i$.

\vspace{0.5em}

\noindent
\textbf{Case 2.} Suppose that the event $\cB_i(v)$ holds.
Here, the new information revealed is the outcomes of the events $\deg_{G_j}(v,V_j(v,i)) < L$ and the neighbours of $v$ in $G(\tau_k)[V \setminus B]$.
As in the previous case, the former events remain independent conditionally on $\cF_{i-1}$. 
Moreover, revealing the indicators of the edges incident to the vertex $v$ and deciding the extensions of the matchings $M_1,\ldots,M_k$ (which is done in a deterministic way conditionally on the said edges) acts as (coupled) one-step exploration processes in each of the graphs $G_1,\ldots,G_k$ and, as such, also preserves the conditional independence of the remaining edges. 
In particular, since the random graphs $(G_j)_{j=1}^k$ were jointly independent conditionally on $\cF_{i-1}$, they remain jointly independent conditionally on $\cF_i$.

\vspace{0.5em}

\noindent
All in all, this justifies the independence of $G_1,\ldots,G_k$ conditionally on $\cF_i$ (and, in particular, the independence of $G_1'',\ldots,G_k''$ under the same conditioning), and concludes the induction.
The analogous independence statement for the matchings follows from the independence of the $G_j''$ and the permutations $\sigma_j$.

\vspace{1em}

Now, suppose that Stage 1 of the matching construction terminates successfully with $V_1,\ldots,V_k,B$.
By definition, none of these sets has changed over the last $n$ steps, which implies that, for every $j\in [k]$ and for every vertex $v\in V\setminus B$, we have $V_j(v,\tau_{\rm mc})=V_j$.
Hence, for every $j\in [k]$, the distribution of $G_j$ conditionally on $\cF_{\tau_{\rm mc}}$ is the one of a binomial random graph with parameters $n$ and $p_1$ with the following restrictions:
\begin{itemize}
\item the set of edges and non-edges incident to the set $B$ is entirely exposed, and
\item for all remaining vertices $v\in V\setminus B$ we have that $\deg_{G_j}(v,V_j) \ge L$.
\end{itemize}
In particular, this distribution is invariant under permutation of the vertices of $V_j$ and thus the distribution of $G_j''$ conditioned on $\cF_{\tau_{\rm mc}}$ is invariant under $\sigma_j$.
This justifies our second claim for $G_1'',\ldots,G_k''$; the uniformity of the matchings now follows directly due to the uniformly random relabelling ensured by $\sigma_j$.
\end{proof}

Next, we show that the construction terminates successfully.
Denote by $S \subseteq V$ the set of vertices with degree at least $2L=2\eta \log n$ in each of the graphs $G_1,\ldots,G_k$ (in the beginning of the process), and set $U = V \setminus S$.
Define an auxiliary graph $\mathcal{A}$ on the vertex set $U$ by adding an edge between any two vertices $u_1, u_2 \in U$ if the shortest path between them in $G(n,q)$ has length at most $4$. 

\begin{prop}\label{prop:A}
Whp every connected component in $\mathcal{A}$ has order at most $15/\delta_1$.
\end{prop}
\proof 
Suppose that the auxiliary graph contains a connected set of $t:=15/\delta$ vertices and note that this implies the existence of a set $W$ which induces a connected graph in $G(n,q)$, has size $w$ with $w \le t+3(t-1) \le 4t$ and contains a subset $W' \subseteq W$ of size $t$ of vertices of $U$; denote the number of the latter sets by $X$. 
Observe that for each $v \in W'$, by definition of $U$, there must exist $i \in [k]$ such that $\deg_{G_i}(v,V) \le 2L$. 
However the events $\{W \text{ is connected in } G(n,q)\}$ and $\{\exists i \in [k]:\deg_{G_i}(v,V) <2L\}_{v \in W'}$ are not independent as the $G_i$ and $G(n,q)$ are coupled with $G(\tau_k)$.
Therefore, instead of $\deg_{G_i}(v,V)$, we consider the random variable $\deg_{G_i}(v,V \setminus W)$, noticing that this is independent of $W$ being connected in $G(n,q)$ and that if $\deg_{G_i}(v,V)< 2L$ then $\deg_{G_i}(v,V\setminus W)< 2L$.
Combining this with the fact that there are $w^{w-2}$ trees on $w$ vertices by Cayley's formula and fixing sets $W'\subseteq W$ with $|W'|=t$ and $|W|=w$, we get
\begin{equation}\label{eq:ubw}
\mathbb E[X] \le \sum_{w=t}^{4t} \binom{n}{w} w^{w-2} q^{w-1} \binom{w}{t} \prob{\forall v\in W',\, \exists i\in [k],\, \deg_{G_i}(v,V \setminus W)< 2L}.
\end{equation}
For the latter probability note that, first, the events $\{\exists i\in [k],\, \deg_{G_i}(v,V\setminus W)< 2L\}$ are jointly independent for the $t$ vertices in $W'$ and, second, for every vertex $v\in V$ and for every graph $G_i$, the degree of $v$ in $V\setminus W$ is distributed as $\mathrm{Bin}(n-w,p_1)$. Since $p_1 = \delta_1 (\log n)/n$ and $\eta \ll \delta_1$, the Chernoff bound and a union bound imply that for any $v \in W'$
\begin{equation}\label{eq:CbUb}
\prob{\exists i\in [k],\, \deg_{G_i}(v,V\setminus W)< 2L}\le 2k\exp\left( - \frac{(\delta_1 \log n-o(1)-2L)^2}{2(\delta_1 \log n + \delta_1 \log n/3)} \right) \le n^{-\delta_1/3}.
\end{equation}

By the independence, $\prob{\forall v\in W',\, \exists i\in [k],\, \deg_{G_i}(v,V \setminus W)\le 2L} \le n^{-t \delta_1/3}=n^{-5}$.
Combined with~\eqref{eq:ubw}, we derive that
\begin{align*}
    \expn{X} \le \sum_{w=t}^{4t} O_t\left(n^w q^{w-1} n^{-5}\right)
    = O_t(n (nq)^{4t-1} n^{-5}) \le n^{-3}\, ,
\end{align*}
where the latter inequality holds for $n$ large enough.
In particular $\expn{X} \to 0$ and thus, by Markov's inequality, whp $\cA$ does not contain a connected component of size $15/\delta_1$, as desired. 
\endproof

We now show that the construction terminates successfully with a small set of bad vertices.

\begin{lemma}\label{lem:B}
Whp the multi-core matching construction terminates successfully. 
Moreover, whp $B\subseteq U$ and the construction ends with $|B|\le |U|\le n^{1-\eps}$.
\end{lemma}
\begin{proof}
Denote by $\cC$ the event from Proposition~\ref{prop:A}. 
We show that whp $B \subseteq U$, meaning that no vertex in $S$ ever enters the set $B$ during the first stage of the multi-core matching construction. 
Suppose for a contradiction that some vertex in $S$ becomes bad, and let $v \in S$ be the \emph{first} such vertex.
Then, the degree of $v$ in some of the graphs $G_i[V \setminus V(M_i)]$ must drop below $L$. Since $v \in S$, it initially possessed at least $2L$ neighbours in $G_i$, meaning it must have lost strictly more than $2L-L = L$ neighbours. 

A neighbour $z$ of $v$ in $G_i$ is lost only when it is added to $V(M_i)$, which occurs exclusively if $z$ itself becomes bad or is chosen as a matching partner in some of $M_1,\ldots,M_k$ of some bad vertex~$z'$.
Because $v$ is chosen as the first vertex in $S$ to turn bad, every bad vertex causing $v$ to lose one neighbour (so the vertex $z$ in the former case and $z'$ in the latter) must belong to $U$.
Moreover, in the first case, $v$ and $z$ are at distance $1$ in $G(n,q)$ and, in the second case, $v$ and $z'$ are at distance at most $2$ in $G(n,q)$.
In particular, the distance between any two such vertices in $G(n,q)$ is at most $4$ and thus they all belong to the same connected component of $\cA$. 
Consequently, unless the event $\cC$ fails (which happens with probability $o(1)$ by Proposition~\ref{prop:A}), $v$ may have lost at most $15/\delta_1 < L$ neighbours in the process, contradicting that it has become bad.

In order to show that $|B| \le n^{1-\eps}$, it is therefore sufficient to show that $|U| \le n^{1-\eps}$. 
By another Chernoff bound and a union bound (applied as in~\eqref{eq:CbUb}), we have $\expn{|U|}\le n^{1-\delta_1/3}$.
Choosing $\eps < \delta_1/3$ and applying Markov's inequality to $|U|$ implies that whp $|U|\le n^{1-\eps}$.

\vspace{0.5em}

Next, we ensure that the first stage of the multi-core matching construction terminates successfully whp. 
Set $\xi:=(k+1) \cdot 15/\delta_1$.
Let $v$ be the bad vertex at which the process terminates unsuccessfully.
We first argue that $v$ has degree at most $\xi$ in $H$.
Indeed, by repeating the argument above, unless the event $\cC$ fails, the vertex $v$ loses at most $15/\delta_1$ neighbours in each $G_i$ for $i \in [k]$ over the process.
If $v$ has degree more than $\xi$ in $H$ (and thus in $G(\tau_k)$) then, at the point $v$ becomes bad, it has still $\xi - k\cdot 15/\delta_1 = 15/\delta_1 > k$ available neighbours in $G(\tau_k)$ which we can use to extend the matchings which do not yet saturate $v$ (and saturate $v$ in each of them).
Therefore the process does not terminate unsuccessfully because of a vertex of degree larger than $\xi$ in $H$.

In order to prove that the process does not terminate unsuccessfully at a vertex of degree at most $\xi$ in $H$, we are going to show that whp for every vertex $v$ with degree at most $\xi$ in $H$ and every vertex $u$ at distance at most $2$ in $G(n,q)$ from $v$, it holds that $u \in S$.
This indeed implies that $v$ does not lose any neighbours and the process does not terminate unsuccessfully at $v$.
We are left to show the former part.
Observe that for fixed distinct vertices $u,w,v$, the events 
\[\big\{uw,wv \in E(G(n,q))\big\},\big\{vs \in E(H) \text{ for } s \in V \setminus \{u,w,v\}\big\}\;\text{and}\;\big\{us \in E(G_i) \text{ for } s \in V \setminus \{u,w,v\}\big\}_{i \in [k]}\] 
are mutually independent.
Moreover, if $\deg_H(v) \le \xi$ then $\deg_H(v,V \setminus \{u,w\}) \le \xi$ and, similarly, if $u \not\in S$, then $\deg_{G_i}(u) < 2L$ for some $i \in [k]$ and thus $\deg_{G_i}(u,V\setminus \{v,w\}) < 2L$.
Therefore the probability that $v$ has degree at most $\xi$ in $H$, $u \not\in S$ and $uwv$ is a path of $G(n,q)$ is at most
\begin{align*}
q^2 &\cdot \mathbb P[\mathrm{Bin}(n-3,p_H)\le \xi]\cdot k \prob{\mathrm{Bin}(n-3,p_1) < 2L}\, .
\end{align*}
By a union bound, the probability that there exist two vertices $v,u$ at distance at most $2$ in $G(n,q)$ and such that $v$ has degree at most $\xi$ in $H$ and $u \not\in S$ is at most
\begin{align*}
\sum_{i=1}^2 &\binom{n}{i+1} q^i\cdot \mathbb P[\mathrm{Bin}(n-3,p_H)\le \xi]\cdot k \prob{\mathrm{Bin}(n-3,p_1) < 2L}\\
&= 2n (nq)^2 \cdot (\xi+1) (n p_H)^\xi (1 - p_H)^{n-3-\xi} \cdot 2ke^{-np_1/3} = o(1) \, ,
\end{align*}
which can be justified as follows.
For $mr >1$, we have 
\[\prob{\mathrm{Bin}(m,r) \le a} = \sum_{a'=0}^a \prob{\mathrm{Bin}(m,r) = a'} \le (a+1)(mr)^a (1-r)^{m-a},\] 
which we used to estimate the first probability applied with $m=n-3$, $r=p_H$ and $a=\xi$, while the second probability has been estimated with a Chernoff's bound.
Then we used that $p_H = 1-(1-p_1)^k \le kp_1$ (which follows from Bernoulli's inequality) and $1-p_H=(1-p_1)^k$, from which we have
\[(nq)^2 \cdot (n p_H)^{\xi} = O((nq)^2 (n p_1)^\xi) = O((\log n)^{\xi+2})\]
and
\[ (1-p_H)^{n-3-\xi} e^{-np_1/3} = (1-p_1)^{(n-3-\xi)k} e^{-np_1/3} \le e^{-np_1(k+1/3)-(1-\xi)k} = n^{-\frac{k+1/3}{k+1/4}+o(1)}.\]
This shows that the first stage of the multi-core matching construction terminates successfully whp.
The second one begins with subgraphs $G_1'' = G_1[V\setminus V(M_1)],\ldots,G_k'' = G_k[V\setminus V(M_k)]$ and every vertex in $V\setminus B$ has at least $L=\eta \log n$ neighbours in each of the graphs $G_1'',\ldots,G_k''$ it belongs to.
By Proposition~\ref{prop:min_degree_perf_matching}, whp they all have a perfect matching and thus the second stage also terminates successfully, which finishes the proof.
\end{proof}

Lemma~\ref{lem:unif} and Lemma~\ref{lem:B} already allow to prove Theorem~\ref{thm:coupling_main}. 
We finish the subsection with more useful observations.
The collection of matchings obtained from the multi-core construction does not yet satisfy the Kotzig property.
Indeed, the union of any two of them may consist of many cycles;
in fact, they can even share the same edge, which we identify with a cycle of length~$2$.
In Section~\ref{sec:merging}, we will replace some of their edges with edges of the random reservoir $G_{k+1}[V \setminus B]$ to merge the various cycles into a Hamilton cycle.
For that, we want to control how the cycles in the union behave with respect to the set $B$.

\begin{lemma}\label{lem:properties}
Let $M_1,\dots,M_k$ be the matchings obtained after a successful termination of the multi-core construction.
Whp, for all pairs of distinct $i,j\in [k]$, each of the following properties holds:
\begin{enumerate}[label=\upshape{\textbf{P\arabic*}}]
\item\label{item:P1} no cycle in $M_i\cup M_j$ of length at most $\sqrt{\log n}$ contains a vertex in $B$,
\item\label{item:P2} every cycle in $M_i\cup M_j$ of length at least $\sqrt{\log n}$ contains three consecutive vertices outside $B$,
\item\label{item:P3} $M_i\cup M_j$ contains 
at most $(\log n)^2$ cycles in total.
\end{enumerate}
\end{lemma}

\proof 
We fix distinct $i, j \in [k]$ and denote by $\cD$ the event that all properties from the statement of Proposition~\ref{prop:A} and Lemma~\ref{lem:B} are satisfied.
Denote by $M_i'\subseteq M_i$ and $M_j'\subseteq M_j$ the partial matchings constructed during the first stage, and define $K = M_i' \cup M_j'$.
To prove \ref{item:P1}, we show that no cycle in $M_i \cup M_j$ of length at most $\sqrt{\log n}$ contains a vertex in $U$; note that this is enough to conclude unless the event $\cD$ fails (which happens with probability $o(1)$).
To begin with, cycles of length $2$ with at least one vertex in $B$ do not appear since all edges in $M_1,\ldots,M_k$ incident to $B$ are distinct by construction.
For longer cycles, 
first observe that for fixed vertices $u_1,\dots,u_\ell$, the event 
\[\big\{u_1 \dots u_\ell u_1 \text{ is a cycle of } G(n,q)\big\}\text{ is independent from }\big\{ u_1s \in E(G_i) \text{ for } s \in V \setminus \{u_\ell, u_1, u_2\}\big\}_{i \in [k]},\] 
and that if $\deg_{G_i}(u_1) < 2L$ for some $i \in [k]$ then $\deg_{G_i}(u_1,V \setminus \{u_\ell, u_1, u_2\}) <2L$.
Therefore the expected number of cycles in $G(n,q)\supset M_i\cup M_j$ of length $\ell\in [3,\sqrt{\log n}]$ which contain a vertex with degree smaller than $2L$ in some of $G_1,\ldots,G_k$ is at most
\[\sum_{\ell=3}^{\sqrt{\log n}} n^\ell q^\ell \cdot \ell \cdot k \prob{\mathrm{Bin}(n-3,p_1)< 2L}\le \sqrt{\log n}\cdot (nq)^{\sqrt{\log n}}\cdot k 2n^{-\delta/3} = o(1)\, ,\]
where we used the Chernoff bound to estimate the probability.
Thus, \ref{item:P1} follows from the previous bound and Markov's inequality.

We now show \ref{item:P2} by starting from an arbitrary vertex and exploring its cycle in $M_i\cup M_j$, which we can assume having length at least $\sqrt{\log n}$ (otherwise, there is nothing to prove).
If we find three consecutive vertices not in $V(K)$, then these vertices are also not in $B$ and we are done.
So, in order to violate \ref{item:P2}, every time the exploration of the cycle leaves $V(K)$, it must return to $V(K)$ within three steps.
On the other hand, if three consecutive vertices belong to $V(K)$, at least one of them must belong to $B$.
Therefore, in any sequence of consecutive vertices belonging to $V(K)$, those in $B$ (at least a third of them) belong to the same component of $\cA$.
Also, unless $\cD$ fails, the length of any such sequence is at most $3\cdot 15/\delta_1=45/\delta_1$.
By Lemma~\ref{lem:unif}, conditioned on Stage $1$ of the process, we have that $M_i''$ is a uniformly random matching.
Therefore, unless $\cD$ fails, each time the exploration leaves $V(H)$, the probability to return to $V(K)$ at a single step is at most $|U|/(n-\sqrt{\log n})\le n^{-\eps/2}/3$: in particular, the probability to return to $V(K)$ within three consecutive steps in this interval is at most $n^{-\eps/2}$.
As a result, the probability to violate \ref{item:P2} within the first $\sqrt{\log n}$ steps for the chosen starting vertex is at most
\[(n^{-\eps/2})^{\sqrt{\log n}/(45/\delta_1 + 3)} = n^{-\omega(1)},\]
where we used that we can find at least $\sqrt{\log n}/(45/\delta_1 + 3)$ pairwise vertex-disjoint subpaths each of length $45/\delta_1 + 3$ in the first $\sqrt{\log n}$ steps and we can apply the argument above on each of them.
A union bound over all $n$ choices for a starting vertex shows that \ref{item:P2} indeed holds whp.

We turn to \ref{item:P3}. The main idea here is to erase cycles and to contract (maximal) paths of even length in $K$ to single vertices first; this identification allows us to bypass the vertices saturated in the first stage of the multi-core matching procedure and work with independent perfect matchings on the reduced vertex set $V'$ (which we, somewhat abusively, continue to denote $M_i$ and $M_j$).
Note that, in the case of a path of odd length $2\ell+1$, we contract an arbitrary subpath of length $2\ell$.

Recall that, unless the event $\cD$ fails, $|V(K)| \le 3n^{1-\eps}$ and thus $|V'|\ge n - 3n^{1-\eps}$. Set $n':=|V'|/2$. We first focus on estimating the cycle count via another exploration procedure.
Namely, fix a uniformly random vertex $v$ in $V'$ and explore its cycle in $M_i\cup M_j$ one vertex at a time, starting with the edge in $M_i$ incident to $v$.
When the cycle of $v$ is closed, check if there are unexplored vertices remaining in $V'$ and, if so, repeat the procedure.
Then, the probability to close a cycle at odd steps is 0 and, by Lemma~\ref{lem:unif}, the probability to close a cycle at the $2i$-th step is $1/(2n'-2i+1)$ for every $i\in [n']$, independently of previous steps.
Denote by $I_i$ the indicator random variable of the latter event.
Then, by combining the fact that $\expn{I_1+\ldots+I_{n'}} = \Theta(\log n)$ with Markov's inequality, 
\[\prob{I_1+\ldots+I_{n'}\ge (\log n)^2}=o(1)\, ,\]
which implies \ref{item:P3}.
\endproof

\subsection{The matchings are shuffled}

We now introduce a definition reflecting the idea that every matching among those constructed via the multi-core matching construction is `shuffled' approximately uniformly and independently of each of the other matchings.
For more intuition we refer to the proof outline in Section~\ref{subsec:outline}.
We say that an oriented path with endpoints $s$ and $t$ is \defn{consistently oriented} if it is oriented from $s$ to $t$ or from $t$ to $s$.

\begin{defin}[$\alpha$-shuffled matchings]
\label{def:shuffled}
    For $\alpha >0$ and a collection of $k$ perfect matchings $M_1, \dots, M_k$ on the same vertex set, we say that the collection is \defn{$\alpha$-shuffled} if the following holds for every $i \in [k]$:
    For every $j \in [k] \setminus \{i\}$, fix a path $P_j$ of $M_i \cup M_j$ with at least $n^{1-\alpha}$ edges of $M_i$ and a consistent orientation of the edges of $P_j$.
    Then 
    \begin{enumerate}[label = (\ding{69})]
		\item \label{item:new_shuffled}
        $\cap_{j \in [k] \setminus \{i\}} P_j$ contains at least $n^{1-\alpha (k-1)}/4^k$ edges of $M_i$ which have the same orientation in each path $P_j$.
    \end{enumerate}
\end{defin}

Next, we show that the family of matchings constructed in the previous section is typically $\alpha$-shuffled for a suitable choice of $\alpha$.
In fact we state a more general lemma to highlight what properties suffice to show that a family is shuffled. 

\begin{lemma}
\label{lem:statistics_new}
    Let $k \ge 2$, $0 < \eps <1$ and $0 < \alpha < \min\{\eps,1/(2k-2)\}$.
    Let $V$ be a set of $n$ vertices.
    For each $i \in [k]$, let $M_i'$ be a matching of size at most $n^{1-\eps}$ and $M_i''$ be a uniformly random perfect matching on $V \setminus V(M_i')$, independently of any other $j \in [k] \setminus \{i\}$.
    Finally let $M_i:=M_i' \cup M_i''$.
    Then whp~the family $\{M_1, \dots, M_k\}$ is $\alpha$-shuffled.
\end{lemma}

\begin{proof}
Fix $i=k$; the arguments for other values of $i\in [k]$ follow from symmetric considerations. 
Colour the edges of $M_k$ in \emph{red} and, for every vertex $v$, denote by $e_v \in M_k$ the red edge incident to $v$. 

Fix a tuple of vertices $(v_1,\dots,v_{k-1})$ and assume that for each $j$ the longest alternating path of $M_i \cup M_j$ starting at $v_j$ with a red edge contains at least $n^{1-\alpha}$ red edges.
We show that with probability $1-o(n^{1-k})$ condition~\ref{item:new_shuffled} holds with $P_j$ being the alternating path of $M_i \cup M_j$ starting at $v_j$ with a red edge, ending with a red edge, having exactly $\lceil n^{1-\alpha}\rceil$ red edges and oriented away from $v_j$.
Then a union bound over all $O(n^{k-1})$ choices of tuples implies the lemma as any consistently oriented path appearing in Definition~\ref{def:shuffled} contains an oriented path which can be constructed as above.

Expose $M_k''$ so that the red matching $M_k$ is fully revealed and, for each colour $j \in [k-1]$, denote by $R_j \subseteq M_k$ the set containing $e_{v_j}$ and all the red edges with at least one endpoint in $V(M_j')$. 
Then, the subgraph $R_j \cup M_j'$ consists of alternating paths and cycles, which we call \defn{$j$-fixed}. 
Observe that each $j$-fixed path starts and ends with a red edge (see~the rectangle in Figure~\ref{fig:permuation} for an example).
Moreover, we call a red edge \defn{$j$-free} if it belongs to $M_k \setminus R_j$, and \emph{free} if it is $j$-free for every $j\in [k-1]$.
(We note that $e_{v_{j}}$ is not $j$-free.)
Observe that $|R_j| \le |V(M_j')|+1 \le 2n^{1-\eps}$ and thus the number of $j$-free edges is at least $n/2-O(n^{1-\eps})$.
It follows that if $t$ denotes the number of free edges, we have $t=n/2-O(kn^{1-\eps})=n/2-o(n^{1-\alpha})$.

Next, for every $j\in [k-1]$, given $M_j'\cup M_k$, we consecutively reveal the edges of $M_j''$ as follows and, simultaneously, we orient the edges of $M_k$.
Observe that $M_j' \cup M_k$ is a collection of alternating cycles and paths, with paths starting and ending with a red edge and possibly being just a single red edge.
We start by orienting the path or cycle containing $e_{v_j}$ away from $v_j$.
As long as the currently explored component of $M_j\cup M_k$ is a path and not a cycle, consider its last vertex $v$; note that the last edge of the path is red and oriented towards $v$.
Next, sample the partner $u$ of $v$ in $M_j''$. 
Two scenarios are possible:
\begin{itemize}
\item If the edge $uv$ closes a cycle, orient the edge from $v$ to $u$.
Then find an arbitrary path $P$ in $M_j'\cup M_k$ and an arbitrary endpoint $w$ of $P$ whose partner in $M_j''$ has not yet been revealed (note that $P$ may be a single $j$-free red edge). If such does not exist, terminate the procedure and output $M_j$. Otherwise, orient $P$ towards $w$.
\item If the edge $uv$ extends the current path (even just by a single $j$-free red edge), then 
orient the newly added path from $u$ to its other endpoint $w$.
\end{itemize}
In each case, redo the procedure with $w$ replacing $v$.
By the end of the procedure, only the $j$-fixed cycles of $M_j' \cup M_k$ remain unexplored and we orient them uniformly at random and independently of the past.
Observe that every red edge is oriented and call this the \emph{$j$-orientation} of $M_j'\cup M_k$.
Moreover note that the free red edges are oriented uniformly at random and independently for each $j \in [k-1]$.

Label the $t$ free red edges arbitrarily and bijectively by integers in $[t]$.
Consider the $k-1$ permutations $\sigma_1,\ldots,\sigma_{k-1}\in S_t$ where $\sigma_j$ indicates the order in which the free edges appear in the exploration procedure of $M_j''$ (see~Figure~\ref{fig:permuation} for an example).
Note that they are independent permutations of $S_t$ and also independent of the $j$-orientations of the free edges for all $j \in [k-1]$.

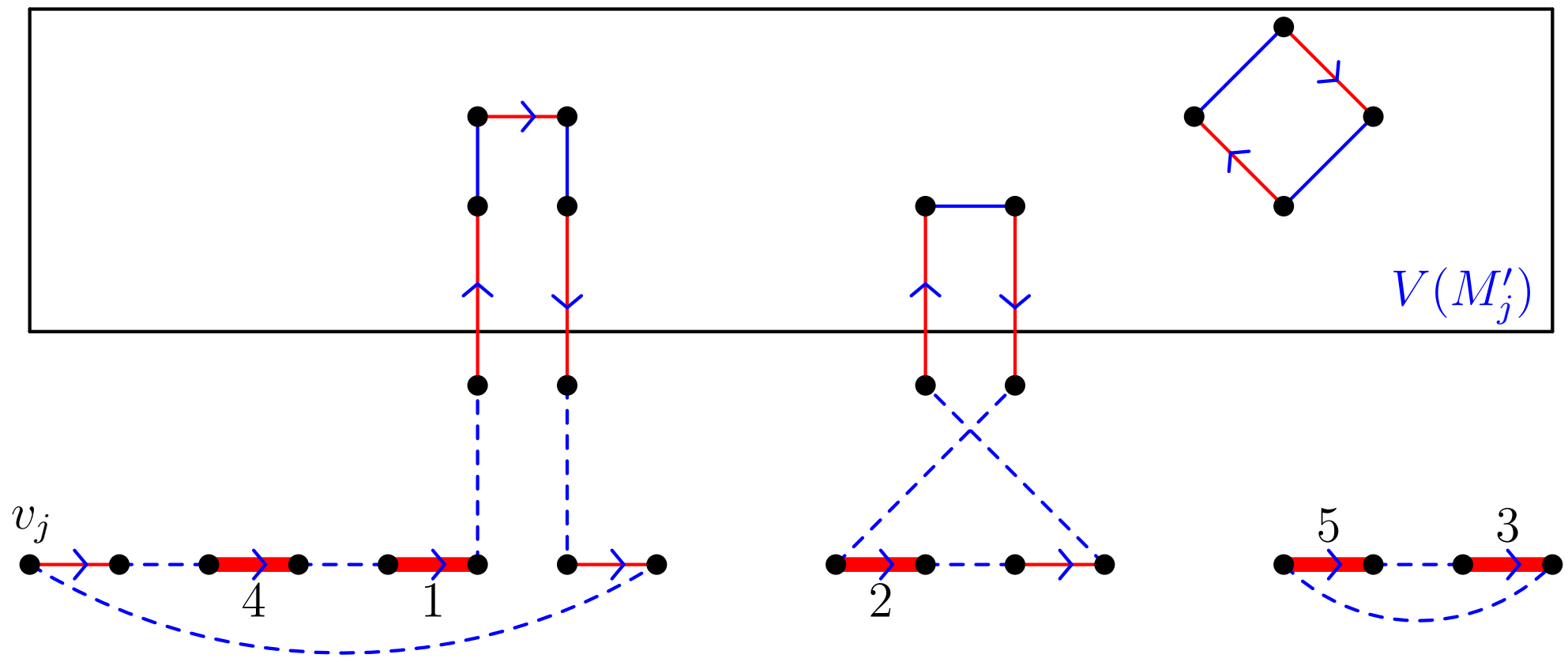
\begin{figure}[htp]
    \centering
\begin{tikzpicture}[scale=1.5,line cap=round,line join=round,x=1cm,y=1cm]
\clip(-7,-0.1) rectangle (2.5,3.7);

\draw [line width=3.5pt,color=red] (-5.5,0.5)-- (-5,0.5); % 4
\draw [line width=3.5pt,color=red] (-4.5,0.5)-- (-4,0.5); % 1
\draw [line width=3.5pt,color=red] (-2,0.5)-- (-1.5,0.5); % 2
\draw [line width=3.5pt,color=red] (0.5,0.5)-- (1,0.5);   % 5
\draw [line width=3.5pt,color=red] (1.5,0.5)-- (2,0.5);   % 3

\draw [line width=0.8pt,color=black] (-6.5,3.6)-- (-6.5,1.8);
\draw [line width=0.8pt,color=black] (-6.5,1.8)-- (2,1.8);
\draw [line width=0.8pt,color=black] (2,1.8)-- (2,3.6);
\draw [line width=0.8pt,color=black] (2,3.6)-- (-6.5,3.6);

\draw [line width=0.8pt,color=red] (0.5,2.5)-- (0,3);
\draw [line width=0.8pt,color=blue] (0.20322632625401876,2.7967736737459807) -- (0.3061284084951766,2.8061284084951774);
\draw [line width=0.8pt,color=blue] (0.20322632625401876,2.7967736737459807) -- (0.19387159150482297,2.693871591504823);
\draw [line width=0.8pt,color=blue] (0,3)-- (0.5,3.5);

\draw [line width=0.8pt,color=red] (0.5,3.5)-- (1,3);
\draw [line width=0.8pt,color=blue] (0.7967736737459803,3.203226326254019) -- (0.6938715915048225,3.1938715915048226);
\draw [line width=0.8pt,color=blue] (0.7967736737459803,3.203226326254019) -- (0.8061284084951761,3.3061284084951765);
\draw [line width=0.8pt,color=blue] (1,3)-- (0.5,2.5);

\draw [line width=0.8pt,color=red] (-4,2.5)-- (-4,1.5);
\draw [line width=0.8pt,color=blue] (-4,2.0661479637735805) -- (-4.079377556528297,2);
\draw [line width=0.8pt,color=blue] (-4,2.0661479637735805) -- (-3.920622443471703,2);
\draw [line width=0.8pt,color=blue] (-4,2.5)-- (-4,3);

\draw [line width=0.8pt,color=red] (-4,3)-- (-3.5,3);
\draw [line width=0.8pt,color=blue] (-3.6838520362264195,3) -- (-3.75,2.920622443471703);
\draw [line width=0.8pt,color=blue] (-3.6838520362264195,3) -- (-3.75,3.0793775565282964);
\draw [line width=0.8pt,color=blue] (-3.5,3)-- (-3.5,2.5);

\draw [line width=0.8pt,color=red] (-3.5,2.5)-- (-3.5,1.5);
\draw [line width=0.8pt,color=blue] (-3.5,1.9338520362264195) -- (-3.579377556528297,2);
\draw [line width=0.8pt,color=blue] (-3.5,1.9338520362264195) -- (-3.4206224434717027,2);

\draw [line width=0.8pt,color=red] (-1.5,1.5)-- (-1.5,2.5);
\draw [line width=0.8pt,color=blue] (-1.5,2.0661479637735805) -- (-1.4206224434717025,2);
\draw [line width=0.8pt,color=blue] (-1.5,2.0661479637735805) -- (-1.5793775565282968,2);
\draw [line width=0.8pt,color=blue] (-1.5,2.5)-- (-1,2.5);

\draw [line width=0.8pt,color=red] (-1,2.5)-- (-1,1.5);
\draw [line width=0.8pt,color=blue] (-1,1.9338520362264195) -- (-1.0793775565282973,2);
\draw [line width=0.8pt,color=blue] (-1,1.9338520362264195) -- (-0.920622443471703,2);

\draw [line width=0.8pt,color=red] (-4.5,0.5)-- (-4,0.5);
\draw [line width=0.8pt,color=blue] (-4.183852036226419,0.5) -- (-4.25,0.42062244347170313);
\draw [line width=0.8pt,color=blue] (-4.183852036226419,0.5) -- (-4.25,0.5793775565282966);

\draw [line width=0.8pt,color=red] (-5.5,0.5)-- (-5,0.5);
\draw [line width=0.8pt,color=blue] (-5.183852036226418,0.5) -- (-5.25,0.42062244347170313);
\draw [line width=0.8pt,color=blue] (-5.183852036226418,0.5) -- (-5.25,0.5793775565282966);

\draw [line width=0.8pt,color=red] (-3.5,0.5)-- (-3,0.5);
\draw [line width=0.8pt,color=blue] (-3.183852036226418,0.5) -- (-3.25,0.42062244347170313);
\draw [line width=0.8pt,color=blue] (-3.183852036226418,0.5) -- (-3.25,0.5793775565282966);

\draw [line width=0.8pt,color=red] (-6.5,0.5)-- (-6,0.5);
\draw [line width=0.8pt,color=blue] (-6.183852036226419,0.5) -- (-6.25,0.42062244347170313);
\draw [line width=0.8pt,color=blue] (-6.183852036226419,0.5) -- (-6.25,0.5793775565282966);

\draw [line width=0.8pt,color=red] (-2,0.5)-- (-1.5,0.5);
\draw [line width=0.8pt,color=blue] (-1.6838520362264195,0.5) -- (-1.75,0.42062244347170313);
\draw [line width=0.8pt,color=blue] (-1.6838520362264195,0.5) -- (-1.75,0.5793775565282966);

\draw [line width=0.8pt,color=red] (-1,0.5)-- (-0.5,0.5);
\draw [line width=0.8pt,color=blue] (-0.6838520362264184,0.5) -- (-0.75,0.42062244347170313);
\draw [line width=0.8pt,color=blue] (-0.6838520362264184,0.5) -- (-0.75,0.5793775565282966);

\draw [line width=0.8pt,color=red] (0.5,0.5)-- (1,0.5);
\draw [line width=0.8pt,color=blue] (0.8161479637735802,0.5) -- (0.75,0.42062244347170313);
\draw [line width=0.8pt,color=blue] (0.8161479637735802,0.5) -- (0.75,0.5793775565282966);

\draw [line width=0.8pt,color=red] (1.5,0.5)-- (2,0.5);
\draw [line width=0.8pt,color=blue] (1.8161479637735813,0.5) -- (1.75,0.42062244347170313);
\draw [line width=0.8pt,color=blue] (1.8161479637735813,0.5) -- (1.75,0.5793775565282966);

\draw [line width=0.8pt,dash pattern=on 3pt off 3pt,color=blue] (-1.5,0.5)-- (-1,0.5);
\draw [line width=0.8pt,dash pattern=on 3pt off 3pt,color=blue] (-0.5,0.5)-- (-1.5,1.5);
\draw [line width=0.8pt,dash pattern=on 3pt off 3pt,color=blue] (-2,0.5)-- (-1,1.5);
\draw [line width=0.8pt,dash pattern=on 3pt off 3pt,color=blue] (-6,0.5)-- (-5.5,0.5);
\draw [line width=0.8pt,dash pattern=on 3pt off 3pt,color=blue] (-5,0.5)-- (-4.5,0.5);
\draw [line width=0.8pt,dash pattern=on 3pt off 3pt,color=blue] (-4,0.5)-- (-4,1.5);
\draw [line width=0.8pt,dash pattern=on 3pt off 3pt,color=blue] (-3.5,1.5)-- (-3.5,0.5);
\draw [line width=0.8pt,dash pattern=on 3pt off 3pt,color=blue] (1,0.5)-- (1.5,0.5);

\draw [shift={(-4.75,3.354131148265822)},line width=0.8pt,dash pattern=on 3pt off 3pt,color=blue]
plot[domain=4.162359107290655:5.262418853478724,variable=\t]
({1*3.3479194451929666*cos(\t r)+0*3.3479194451929666*sin(\t r)},
 {0*3.3479194451929666*cos(\t r)+1*3.3479194451929666*sin(\t r)});

\draw [shift={(1.25,1.2406453617582658)},line width=0.8pt,dash pattern=on 3pt off 3pt,color=blue]
plot[domain=3.920715336795219:5.504062623974161,variable=\t]
({1*1.0540661990093565*cos(\t r)+0*1.0540661990093565*sin(\t r)},
 {0*1.0540661990093565*cos(\t r)+1*1.0540661990093565*sin(\t r)});

\begin{scriptsize}
\draw[color=blue] (1.5,2) node {\large{$V(M_j')$}};
\draw [fill=black] (0.5,2.5) circle (1.5pt);
\draw [fill=black] (0,3) circle (1.5pt);
\draw [fill=black] (0.5,3.5) circle (1.5pt);
\draw [fill=black] (1,3) circle (1.5pt);
\draw [fill=black] (-4,2.5) circle (1.5pt);
\draw [fill=black] (-4,1.5) circle (1.5pt);
\draw [fill=black] (-4,3) circle (1.5pt);
\draw [fill=black] (-3.5,3) circle (1.5pt);
\draw [fill=black] (-3.5,2.5) circle (1.5pt);
\draw [fill=black] (-3.5,1.5) circle (1.5pt);
\draw [fill=black] (-1.5,1.5) circle (1.5pt);
\draw [fill=black] (-1.5,2.5) circle (1.5pt);
\draw [fill=black] (-1,2.5) circle (1.5pt);
\draw [fill=black] (-1,1.5) circle (1.5pt);
\draw [fill=black] (-4.5,0.5) circle (1.5pt);
\draw [fill=black] (-4,0.5) circle (1.5pt);
\draw [fill=black] (-5.5,0.5) circle (1.5pt);
\draw [fill=black] (-5,0.5) circle (1.5pt);
\draw [fill=black] (-3.5,0.5) circle (1.5pt);
\draw [fill=black] (-3,0.5) circle (1.5pt);
\draw [fill=black] (-6.5,0.5) circle (1.5pt);
\draw[color=black] (-6.4870885003248056,0.72195336913361) node {\large{$v_j$}};
\draw [fill=black] (-6,0.5) circle (1.5pt);
\draw [fill=black] (-2,0.5) circle (1.5pt);
\draw [fill=black] (-1.5,0.5) circle (1.5pt);
\draw [fill=black] (-1,0.5) circle (1.5pt);
\draw [fill=black] (-0.5,0.5) circle (1.5pt);
\draw [fill=black] (0.5,0.5) circle (1.5pt);
\draw [fill=black] (1,0.5) circle (1.5pt);
\draw [fill=black] (1.5,0.5) circle (1.5pt);
\draw [fill=black] (2,0.5) circle (1.5pt);
\draw[color=black] (-1.75,0.3) node {\large{$2$}};
\draw[color=black] (-5.25,0.3) node {\large{$4$}};
\draw[color=black] (-4.25,0.3) node {\large{$1$}};
\draw[color=black] (0.75,0.72) node {\large{$5$}};
\draw[color=black] (1.75,0.72) node {\large{$3$}};
\end{scriptsize}
\end{tikzpicture}
    \caption{The edges of $M_k$ are depicted in red, those of $M_j$ in blue.
    The rectangle denotes $V(M_j')$ and contains the $j$-fixed paths and cycles.
    The thick red edges are the free edges which are labeled bijectively by integers in $[t]$. 
    We reveal $M_j''$ as described in the text starting from $e_{v_j}$ and giving an orientation to all the red edges.
    In this case, we have 
    \centerline{$\sigma_j=
    \left(
        \begin{matrix}
          1 & 2 & 3 & 4 & 5 \\
          4 & 1 & 2 & 5 & 3 
        \end{matrix}
    \right)$.}}
    \label{fig:permuation}
\end{figure}

Recall that the aim is to count the number of red edges in $\cap_{j \in [k-1]} P_j$ which have the same $j$-orientation for each $j \in [k-1]$.
We start by counting the number of red edges in the intersection, without taking the orientation into account.
Given that the number of red edges which are not free is $n/2-t=o(n^{1-\alpha})$, among the $n^{1-\alpha}$ red edges of $P_j$, at least $n^{1-\alpha}/2$ are free. 
Therefore the first $n^{1-\alpha}/2$ free red edges encountered in the exploration of $M_j''$ belong to $P_j$ and the set of their labels can be expressed as $\{\sigma_j(m):m \in [n^{1-\alpha}/2]\}$.
Now, define the function $f: S_t^{k-1} \to \mathbb{N}$ mapping any tuple of permutations $(\pi_1,\ldots,\pi_{k-1})$ to the size of the set
\[ F(\pi_1,\ldots,\pi_{k-1}) \coloneqq \bigcap_{j=1}^{k-1} \bigg\{ \pi_j(m) : m \in \left[ n^{1-\alpha}/2 \right] \bigg\}. \]
By the argument above, the number of red edges in $\cap_{j \in [k-1]} P_j$ is at least $f(\sigma_1,\dots,\sigma_{k-1})$.
Further note that, for all $j\in [k-1]$ and all permutations $\pi_1,\ldots,\pi_{k-1},\pi_j'$ where  $\pi_j$ and $\pi_j'$ differ by a single transposition,
\[|f(\pi_1,\ldots,\pi_j,\ldots,\pi_{k-1})-f(\pi_1,\ldots,\pi_j',\ldots,\pi_{k-1})|\le 1.\]
Since each permutation is uniformly random, each free red edge has the same probability of appearing in the first $n^{1-\alpha}/2$ free red edges of the exploration and their number is $t$. 
By independence, we also have
\[\expn{f(\sigma_1,\ldots,\sigma_{k-1})} = t\cdot \left(\frac{n^{1-\alpha}/2}{t}\right)^{k-1} \ge n^{1-(k-1)\alpha}/4,\]
where we used that $n/4 \le t \le n/2$.
McDiarmid's inequality for permutations (Lemma~\ref{lem:mcdiarmid_permutation}) implies that
\[\prob{f(\sigma_1,\ldots,\sigma_{k-1})\le n^{1-(k-1)\alpha}/8}\le 2\exp\bigg(-\frac{n^{2-2(k-1)\alpha}/64}{2t(k-1)}\bigg) = o(n^{1-k}),\]
where in the last inequality we used that $\alpha < 1/(2k-2)$.

Now, given any realisation of $\sigma_1,\ldots,\sigma_{k-1}$ satisfying that $|F|\ge n^{1-(k-1)\alpha}/8$, we reveal the random orientations of the free edges in $F$.
Using that every free red edge in $F$ has the same $j$-orientation for every $j\in [k-1]$ with probability $2 \cdot (1/2)^{k-1} = 1/2^{k-2}$, the Chernoff bound implies that, with probability $1-o(n^{1-k})$, there are at least $|F|/2^{k-1}\ge n^{1-(k-1)\alpha}/4^k$ edges in $F$ with the same $j$-orientation for each $j$.
As observed before, these edges belong to $\cap_{j \in [k-1]} P_j$.
Finally, a union bound over the $kn^{k-1}$ choices of vertices $v_1, \dots,v_{k-1}\in V$ and an index $i\in [k]$ concludes the lemma.
\end{proof}

\section{Merging cycles}
\label{sec:merging}

Recall that Theorem~\ref{thm:coupling_main} proves that $G(\tau_k)$ contains a collection of matchings $M_1, \dots, M_k$, a set of vertices $B$ and a set of edges $R$ with certain properties.
In this section, we utilise edges from $R$ to iteratively merge the cycles in each $2$-factor $M_i \cup M_j$ into a Hamilton cycle, following the strategy outlined in Section~\ref{subsec:outline}.

We first explain the edge-exchange operations we do to merge cycles in a concrete 2-factor without fragmenting further other 2-factors. 
Suppose we are given $k$ perfect matchings $M_1,\ldots,M_k$ and we wish to modify $M_i$ so that $M_i \cup M_j$ becomes a Hamilton cycle for each $j \in [k] \setminus \{i\}$.
A graph is called an \defn{$M_i$-near $2$-factor} (for $M_1,\ldots,M_k$) if it is obtained by taking out one edge in $M_i$ from the union of $M_i\cup M_j$ for some $j\in [k]\setminus \{i\}$: in particular, an $M_i$-near 2-factor contains only one path whose endpoints are not saturated by $M_i$ and potentially several cycles alternating between $M_i$ and $M_j$. 
For every $i\in [k]$ and every $M_i$-near $2$-factor $L$, we denote by $P(L)$ its unique path.

Our approach relies on P\'osa-style rotations using edges from the `reservoir' random graph $G_{k+1}$ to extend or alter the $M_i$-near 2-factors. 
However, when replacing edges of $M_i$ with edges of $G_{k+1}$, we must avoid creating more cycles in each of the $M_i$-near 2-factors. The `good' edges supplied by the following definition indicate the places where rotations of the alternating paths can be executed.

\begin{defin}[Good edges]
\label{def:good}
Fix $i\in [k]$, an edge $uv\in M_i$ and consider the $M_i$-near $2$-factors $L_j=M_i\cup M_j\setminus \{uv\}$ for all $j\in [k]\setminus \{i\}$. 
An edge $f=xy\in M_i\setminus \{uv\}$ is called \defn{$(M_i,v)$-good} if, for every $j\in [k]\setminus \{i\}$ such that $f\in P(L_j)$, $f$ is traversed in the same direction when running through $P(L_j)$ from $u$ to $v$.
If this common direction is from $x$ to $y$, we call $x$ the \defn{good endpoint}~of~$f$. 
\end{defin}

We remark that the notion of being a good edge does not depend on $M_i$ and $v$ only, but on the entire family of matchings; we write $(M_i,v)$-good to ease the notation. 
Suppose $xy$ is an $(M_i,v)$-good edge with good endpoint $x$. 
If $vx\in E(G_{k+1})$, we can exchange $xy$ with $vx$ in $M_i$ without increasing the number of components in any $M_i$-near 2-factor: indeed, the number of components remains unchanged for all $M_i \cup M_j$ with $xy\in P(L_j)$ while it decreases by 1 in all other $M_i$-near 2-factors. 
We call this operation an \defn{acceptable rotation around $v$ with fixed endpoint $u$ and pivot $x$}; note that it makes $y$ the new endpoint of the path from $u$ in each $M_i$-near 2-factor.

Next, we show that, for every $\alpha$-shuffled collection of matchings, modifying a polylogarithmic number of edges leaves a large number of acceptable rotations for every $i\in [k]$, even when we require the rotations to use a fixed subsets of edges (which is denoted by $(F_j)_{j\in [k]\setminus \{i\}}$ in the statement of the proposition). 
For a graph $G$, we say that a subgraph $H \subseteq G$ is \defn{component-respecting} if for each component $C$ of $G$ either $C \subseteq H$ or $C \cap H = \emptyset$.

\begin{prop}
\label{prop:stays_shuffled}
Fix $1/n \ll  \alpha, 1/k$ and an $\alpha$-shuffled collection of perfect matchings $M_1,\dots, M_k$ on $V$ such that, for all distinct $i,j \in [k]$, $M_i \cup M_j$ contains at most $(\log n)^2$ cycles.
For each $j \in [k]$, let $\hat M_j$ be a perfect matching on $[n]$ satisfying $|M_j \Delta \hat M_j| \le (\log n)^4$. 
Then, fix $i \in [k]$, an edge $e=uv \in \hat M_i$ and, for each $j \in [k] \setminus \{i\}$, fix a component-respecting subgraph $F_j \subseteq \hat M_i\cup \hat M_j \setminus \{e\}$ with $|V(F_j)| \ge n/2$.
Then, after deleting $e$ from $\hat M_i$, there are at least $n/20^k$ $(\hat M_i,v)$-good edges (for the collection $\hat M_1,\ldots,\hat M_k$) belonging to $\cap_{j \in [k] \setminus \{i\}} F_j$.
\end{prop}

\begin{proof}
We show the statement when $i=k$; the arguments for other values of $i\in [k]$ follow from symmetric considerations. 
Call an edge \emph{red} if it belongs to $M_k\cap \hat M_k$ and, for $j \in [k-1]$, call an edge \emph{$j$-usable} if it belongs to $(M_j \cap \hat M_j) \cup (M_k \cap \hat M_k) \setminus \{e\}$.

Note that the $j$-usable edges form a collection of $m\le |M_j \Delta \hat M_j|+|M_k \Delta \hat M_k|+(\log n)^2\le 3(\log n)^4$ paths and cycles, which cover all but at most $3(\log n)^4$ edges of $\hat M_k \cup \hat M_j$. 
In particular, since $F_j$ is component-respecting, each such path/cycle is either fully contained in $F_j$ or disjoint from $F_j$.
Consequently, the number of $j$-useable edges in $F_j$ is at least $n/2 -3(\log n)^4$, and these edges form a collection of $m' \le m$ paths and cycles.

Thus, for all $j\neq i$, we can greedily extract a collection $\cP_j$ of pairwise edge-disjoint paths of $j$-usable edges in $F_j$, each starting with a red edge and containing exactly $n^{1-\alpha}$ red edges such that $|\cP_j| \ge n^\alpha/5$.
This can be proved as follows.
Let $\ell_1, \dots, \ell_m$ be the lengths of the paths and cycles formed by the $j$-usable edges in $F_j$.
Observe that any cycle or path with $\ell$ edges contains $\lfloor \frac{\ell}{2n^{1-\alpha}}\rfloor$ pairwise edge-disjoint subpaths starting with a red edge and containing exactly $n^{1-\alpha}$ red edges.
Then,
\[|\cP_j| = \sum_{h=1}^{m'} \left\lfloor \frac{\ell_h}{2n^{1-\alpha}}\right\rfloor \ge \sum_{h=1}^{m'} \left(\frac{\ell_h}{2n^{1-\alpha}}-1\right) =  \frac{\sum_{h=1}^{m'} \ell_h}{2n^{1-\alpha}}-m' 
    \ge \frac{n/2 - 3(\log n)^4}{2n^{1-\alpha}}-3(\log n)^4 
    \ge \frac{n^{\alpha}}{5}.\]

Consider any tuple of paths $(P_1, \dots, P_{k-1})$ with $P_j \in \cP_j$ for each $j \in [k-1]$ and observe that each $P_j$ contains $n^{1-\alpha}$ red edges and that they belong to each of $M_k\cup M_j$ and to $F_j$ by definition.
If $P_j$ belongs to the unique path of the $\hat{M_k}$-near $2$-factor $\hat{M_k} \cup \hat{M_j} \setminus \{e\}$, orient $P_j$ from $u$ to $v$; otherwise, orient it arbitrarily.
As the collection $M_1,\ldots,M_k$ is $\alpha$-shuffled, the intersection $\cap_{j \in [k-1]} P_j$ contains at least $n^{1-\alpha(k-1)}/4^k$ red edges oriented the same way for every $j \in [k-1]$. 
These edges belong to $\hat M_k$ and are thus $(\hat M_k,v)$-good; moreover they belong to $\cap_{j \in [k] \setminus \{i\}} F_j$. 
Summing over all $|\cP_1| \cdot \ldots \cdot |\cP_{k-1}|$ choices of tuples, the number of $(\hat M_k,v)$-good edges is at least $(n^{1-\alpha(k-1)}/4^k) \cdot (n^\alpha/5)^{k-1} \ge n/20^k$, as desired.
\end{proof}

To repeatedly apply acceptable rotations, we must ensure that the required pivotal edges exist in the random reservoir $G_{k+1}$. 
We model the sequential exploration of acceptable rotations as a Galton-Watson branching process. 
The next proposition shows a convenient property of this process.

\begin{prop}\label{prop:branching process}
Fix $c_1,c_2 > 0$ and $\gamma\in (0,1)$. Set $p=c_1(\log n)/n$ and consider a branching process $(Z_\ell)_{\ell \ge 0}$ with offspring distribution $\mathrm{Bin}(c_2n, p)$ and initial condition $|Z_0|=1$. 
Then, with probability $1-o((\log n)^{-3})$, there is $h = (\gamma+o(1))\tfrac{\log n}{\log\log n}$ such that 
\[n^{\gamma}\le |Z_h|\le |Z_0\cup Z_1\cup\ldots\cup Z_h|\le n^{\gamma}(\log n)^2.\]
\end{prop}
\begin{proof}
Set $\mu = c_2np$ and denote by $\cE_i$ the event that $|Z_i| = (1\pm 0.5)\mu |Z_{i-1}|$. 
Note that, on the event $\cE_1\cap \ldots\cap \cE_{i-1}$, we have $ (\mu/2)^{i-1} \le |Z_{i-1}| \le (2\mu)^{i-1}$ and, conditionally on $|Z_{i-1}|$, the distribution of $|Z_i|$ is $\mathrm{Bin}(|Z_{i-1}|c_2n,p)$.
By this observation and iterative applications of the Chernoff bound,
\begin{align*}
\prob{\bigcap_{i=1}^{\infty}\cE_i} 
&\ge 1 - \sum_{i=1}^{\infty}\prob{\cE_i^c\mid \cE_1\cap \ldots\cap \cE_{i-1}}\\
&\ge 1-\sum_{i=1}^{\infty}2\exp\bigg(-\frac{(\mu\cdot (\mu/2)^{i-1}/2)^{2}}{2(1+1/6)\mu\cdot (\mu/2)^{i-1}}\bigg)\\
&=1-\sum_{i=1}^{\infty}2\exp\bigg(-\frac{3}{14} \left(\frac{\mu}{2}\right)^i\bigg) =1-o((\log n)^{-3}),
\end{align*}
where we used that $\mu=\Theta(\log n)$.
Condition on the intersection of the events $(\cE_i)_{i=1}^{\infty}$ and 
let $h$ be the smallest index such that $|Z_h| \ge n^\gamma$.
Then $n^\gamma > |Z_{h-1}|\ge (\mu/2)^{h-1}$ and $n^\gamma \le |Z_h| \le (2 \mu)^h$, 
which readily implies that $h = (\gamma+o(1))\tfrac{\log n}{\log\log n}$.
Moreover, $|Z_h| \le 2 \mu|Z_{h-1}| < 2\mu n^\gamma$ and, for $0 \le i \le h$ we have $|Z_i| \le |Z_h| (2/\mu)^{h-i}$.
Therefore $|Z_0 \cup Z_1 \cup \dots \cup Z_{h}| \le |Z_h| \cdot \sum_{i=0}^{h} (2/\mu)^{i} \le 2\mu n^\gamma \cdot \frac{1}{1-2/\mu} \le n^\gamma(\log n)^2$ and $h$ satisfies the condition in the statement.
\end{proof}

Finally we show how an $\alpha$-shuffled family of perfect matchings can be turned into a family satisfying the Kotzig property by using some edges from the random reservoir graph. 

\begin{lemma}
\label{lemma:rotations}
    Let $\eps \ll \delta$ and $\alpha \in (0,1)$.
    Let $V$ be a set of $n$ vertices, $B \subseteq V$ with $|B| \le n^{1-\eps}$ and $\{M_1,\dots,M_k\}$ be an $\alpha$-shuffled family of perfect matchings on $V$ such that the properties \ref{item:P1}--\ref{item:P3} hold.
    Let $p=\delta (\log n)/n$ and $R$ be a random graph on $V \setminus B$ distributed as $G(n-|B|,p)$.
    Then, whp $M_1 \cup \dots \cup M_k \cup R$ contains a family of $k$ perfect matchings with the Kotzig property.
\end{lemma}

\proof
Fix $K:=\tbinom{k}{2}$ and, for all $1 \le i< j \le k$, let $R_{i,j}$ be independent random graphs on $V \setminus B$ distributed as $R_{i,j} \sim G(n-|B|,p/K)$. 
Then, we can couple these random graphs with $R$ so that, with probability one, $\bigcup_{1 \le i<j \le k} R_{i,j} \subseteq R$ since that union is distributed as a random graph on $V \setminus B$ with edge density $1-(1-p/K)^K \le p$.
We merge the cycles in $K$ sequential rounds indexed by the pairs $(i,j)\in [k]^2$ with $i<j$ ordered lexicographically. 
In round $(i,j)$, we modify only the matching $M_i$ using edges from $R_{i,j}$ until the union of $M_i$ and $M_j$ becomes a Hamilton cycle, while making sure not to increase the number of cycles in the union of any pair of other matchings. 
Following our strategy, at any step of the proof, the union of any two matchings will contain at most $(\log n)^2$ cycles and we will have only changed at most $(\log n)^4$ edges of each matching, and thus it will always be possible to apply Proposition~\ref{prop:stays_shuffled}.
In particular, before round $(i,i+1)$, we
will have modified only the matchings $M_{j'}$ with $j' <i$ while the matchings $M_i, M_{i+1},\dots,M_{k}$ have not been changed.

By properties \ref{item:P1} and \ref{item:P2}, there exists a set $X[i,i+1] \subseteq M_i$ which contains an edge of $M_i$ disjoint from $B$ for every cycle of $M_i \cup M_j$ and for every $j \in [i+1,k]$.
Moreover, by \ref{item:P3}, we can guarantee that $|X[i,i+1]| \le k(\log n)^2$.
We assume these two properties in what follows.

Denote by $M_1^-, \dots, M_k^-$ the matchings before round $(i,j)$, and by $M_1^+, \dots, M_k^+$ the matchings after round $(i,j)$.
We preserve the following properties at the end of each round:
    \begin{enumerate}[label={\bfseries Q\arabic*}]
        \item \label{item:Q1} $M_s^+ = M_s^-$ for each $s \in [k] \setminus \{i\}$;
        \item \label{item:Q2} $M_i^+ \subseteq M_i^- \cup R_{i,j}$ and $|M_i^+ \Delta M_i^-|\le (\log n)^3$;
        \item \label{item:Q3} $M_i^+ \cup M_j^+$ is a Hamilton cycle;
        \item \label{item:Q4} for all $s\in [k] \setminus \{i\}$, the vertex set of every cycle in $M_i^+ \cup M_s^+$ is the union of the vertex sets of one or several cycles in $M_i^- \cup M_s^-$;
        \item \label{item:Q5} if $j<k$, there is a set $X[i,j+1]$ of size at most $k (\log n)^2$ which contains an edge of $M_i^+$ disjoint from $B$ for every cycle of $M_i^+ \cup M_{j'}^+$ and for every $j' \in [j+1,k]$.
    \end{enumerate}
    
    In the beginning of round $(i,j)$, observe that the 2-factor $M_i^- \cup M_j^-$ contains $t \le (\log n)^2$ cycles. 
    We merge them into a Hamilton cycle in at most $t-1$ iterations.
    Denote by $t_{s-1}$ the number of cycles in the union of the $i$-th and the $j$-th matching after iteration $s-1$ and note that $t_{s-1} \le t - s + 1$. In particular, if $t_{s-1}=1$, the round $(i,j)$ terminates.
    Assume $t_{s-1} \ge 2$ and let $W_{s-1}$ to be the set of vertices with one of the following properties:
    \begin{enumerate}
    \item[(a)] their entire neighbourhood in $R_{i,j}$ has already been exposed, or
    \item[(b)] their current neighbour in the $i$-th matching is of type (a).
    \end{enumerate}
    These are the vertices which we do not want to use as new endpoints or pivots of any new rotation.
    Moreover we will guarantee that $|W_{s-1}| \le 4(s-1)n^{2/3}(\log n)^2$ and 
    \begin{enumerate}[label={\bfseries Q6}]
        \item \label{item:Q6} we maintain the properties \ref{item:Q1}, \ref{item:Q2} and \ref{item:Q4} after each iteration throughout round $(i,j)$ together with the following stronger version of \ref{item:Q5}, namely that there is a set $X$ of size at most $k(\log n)^2$ which contains an edge of the current $i$-th matching disjoint from $B \cup W_{s-1}$ for each cycle in the union of the current $i$-th and $j'$-th matchings with $j' \in [j,k]$.
    \end{enumerate}
    Observe this is clearly satisfied for $s=1$ by our induction hypothesis.
    For ease of writing, denote by $\hat M_i$ and $\hat M_j$ the two matchings this iteration starts with and by $\hat{X}$ the set guaranteed by \ref{item:Q6}.
    We also maintain a dynamic set $W$ whose original value is $W_{s-1}$ and whose size we maintain smaller or equal than $4sn^{2/3}(\log n)^2$ at all times. Observe that the latter bound is satisfied in the beginning of the iteration.

    Denote by $C$ a shortest cycle in $\hat M_i \cup \hat M_j$ and set $e=uv \in C \cap \hat{X}$.
    For each $j' \in [j+1,k]$, denote by $D_{j'}$ the cycle in $\hat M_i \cup \hat M_{j'}$ containing $e$, if this cycle has length at most $n/2$, and let $D_{j'}=\emptyset$ otherwise.
    As $u, v \notin W \cup B$, we delete $e$ from the $i$-th matching (thus producing a path from $C$) and can start a branching process $(Z_\ell)_{\ell \ge 0}$.
    Starting with $Z_0 = \{v\}$, a vertex $z$ is a child of $y$ in this process if $z \notin B \cup W$ and the following holds.
    Denote by $\hat{M_s}[y]$ the matching obtained by modifying $\hat{M_s}$ by performing the unique sequence of rotations from $v$ to $y$ in the branching process.
    Then, there must be an edge $wz \in \hat M_i[y] \setminus \hat{X}$ such that $wz$ is an $(\hat M_i[y],y)$-good edge (with respect to the family $\{\hat M_1[y], \dots, \hat M_k[y]\}$) with good endpoint $w$ and $yw \in E(R_{i,j})$. For the very first generation ($y=v$), we additionally require that $wz \notin E(C) \cup \big(\bigcup_{j' \in [j+1,k]} E(D_{j'})\big)$.
    (This is a technical assumption whose significance will become clear in the end of the proof.)
    Note that every discovered child $z$ corresponds to a valid acceptable rotation around $y$ with pivot $w$ utilising an edge from $R_{i,j}$
    and that the vertices $z$ and $w$ get added to $W$.

    We bound from below the number of $(\hat M_i[y], y)$-good edges not belonging to $\hat{X}$ and disjoint from $B \cup W$.
    By Proposition~\ref{prop:stays_shuffled} (which we can apply since $|M_i \Delta \hat M_i[y]| \le k(\log n)^3$ by \ref{item:Q2} and \ref{item:Q6}) and the fact that $C$ is component-respecting and $|V(C)| \le n/2$, and $D_{j'}$ is component-respecting and $|V(D_{j'})|\le n/2$ for each $j' \in [j+1,k]$, the number of available $(\hat M_i[y], y)$-good edges in the complement of $C \cup \big(\bigcup_{j' \in [j+1,k]} D_{j'}\big)$ is at least $n/20^k$.
    (Note that the extra requirement that the good edge belongs to that set only concerns the first rotation/extension step at a given iteration.)
    Furthermore, the number of edges belonging to $\hat X$ or incident to $B \cup W$ is at most 
    \[k(\log n)^2 + 2n^{1-\eps} + 4s n^{2/3}(\log n)^2 = o(n).\] 
    Thus, in all cases, the number of valid $(\hat M_i[y],y)$-good endpoints is at least $n/(2\cdot 20^k)$ and we take any subset of them of size exactly $\lfloor n/(2\cdot 20^k) \rfloor$. Because each potential cross-edge $yw$ is present in $R_{i,j}$ independently with probability $p/K$, the offspring distribution of the considered branching process coincides with $\mathrm{Bin}(\lfloor n/(2\cdot 20^k) \rfloor, p/K)$.
    By Proposition~\ref{prop:branching process}, with probability $1-o((\log n)^{-3})$, this branching process reaches some generation $h = (\tfrac{2}{3}+o(1))\tfrac{\log n}{\log \log n}$ where $|Z_h| \ge n^{2/3}$ while exposing edges incident to at most $n^{2/3}(\log n)^2$ vertices. 
    We add the newly exposed vertices in $Z_0\cup\ldots\cup Z_h$ and their (new) neighbours in the $i$-th matching to $W$ so 
    \[|W| \le |W_{s-1}| + 2n^{2/3}(\log n)^2 \le (4(s - 1) + 2)n^{2/3}(\log n)^2.\]
    
    Next, we initiate a second branching process $(Z_{\ell}')_{\ell\ge 0}$ from $u$ defined similarly but with the following modifications. Starting with $Z_0' = \{u\}$, a vertex $z$ is a child of $y$ in this process if $z \notin B \cup W$ and there is an edge $wz \in \hat M_i \setminus \hat X$ such that $yw \in E(R_{i,j})$. 
    As before, $z$ and $w$ get added to $W$.
    By considering only $n/2$ vertices of $V \setminus (B \cup W)$ similarly as above, the offspring distribution of $(Z_{\ell}')_{\ell\ge 0}$ coincides with $\mathrm{Bin}(n/2,p/K)$.
    Note that we no longer require that $wz$ is a good edge with good endpoint $w$. 
    This is because the $\hat M_i$-near $2$-factor obtained after rotations and extensions following the trajectory between $v$ and some $v^*\in Z_h$ in the branching process $(Z_{\ell})_{\ell\ge 0}$ is different for different choices of $v^*$.
    This also influences the sets of good edges, which are now dependent on the endpoint~$v^*$.

    We explain how to resolve the above inconvenience. 
    By Proposition~\ref{prop:branching process}, with probability $1-o((\log n)^{-3})$, there is $h'=(\tfrac{2}{3}+o(1))\tfrac{\log n}{\log \log n}$ with $|Z'_{h'}|\ge n^{2/3}$ and $|Z_0'\cup \ldots\cup Z'_{h'}|\le n^{2/3}(\log n)^2$.
    For every pair of vertices $v^*\in Z_h$ and $u^*\in Z_0'\cup\ldots\cup Z'_{h'}$ and $s \in [k]$, denote by $\hat M_s[v^*,u^*]$ the matching obtained by modifying $\hat M_s$ by performing the unique sequence of rotations from $v$ to $v^*$ in $(Z_{\ell})_{\ell\ge 0}$, and then the unique sequence of rotations from $u$ to $u^*$ in $(Z'_{\ell})_{\ell\ge 0}$.
    This is possible since the branching process $(Z_{\ell}')_{\ell\ge 0}$ stays disjoint from $B\cup W$ and thus it does not depend on the particular trajectory of $(Z_\ell)_{\ell \ge 0}$.
    More precisely, if $z$ is a child of $y$ at any generation of $(Z_{\ell}')_{\ell\ge 0}$ then, after any sequence of rotations in $(Z_\ell)_{\ell \ge 0}$, the edge $wz$ is still an edge of $\hat M_i \setminus \hat X$.

    \begin{cl}\label{cl:good}
    With probability $1-o((\log n)^{-3})$, for at least $n^{2/3}/2$ vertices $v^*\in Z_h$, there are at least $n^{1/2}$ vertices $u^*\in Z'_{h'}$ such that, 
    for every vertex $y$ on the trajectory between $u$ and $u^*$ in $(Z'_{\ell})_{\ell\ge 0}$, before replacing the edge $wz \in \hat M_i[v^*,y]$ with $yw \in E(R_{i,j})$, the edge $wz$ is $(\hat M_i[v^*,y],y)$-good with good endpoint $w$ (with respect to the family $\{\hat M_1[v^*,y], \dots, \hat M_k[v^*,y]\}$).
    \end{cl}
    \lateproof{Claim~\ref{cl:good}}
    By \ref{item:Q2} and \ref{item:Q6}, and the fact that the number of edges where $\hat M_i$ and $\hat M_i[v^*,y]$ differ is at most $h+h' = o(\log n)$, we can apply Proposition~\ref{prop:stays_shuffled}.
    So, for every vertex $y$ on a layer among $0,\ldots,h'-1$ in the branching process $(Z'_{\ell})_{\ell\ge 0}$, there exist at least $n/20^k$ edges which are $(\hat M_i[v^*,y],y)$-good for the current matching configuration.
    Using that the number of edges belonging to $\hat X$ or incident to $W \cup B$ is at most $o(n)$, we deduce that there are at least $n/(2\cdot 20^k)$ good endpoints and thus one can extract a thinner branching process $(Z''_{\ell}[v^*])_{\ell\ge 0}$ from $(Z'_{\ell})_{\ell\ge 0}$ which consists only of acceptable rotations and has
    offspring distribution $\mathrm{Bin}(n/(2 \cdot 20^k),p/K)$.
    Fix $\gamma\in (1/2,2/3)$.
    Since $h' \ge (\gamma+\Omega(1))\frac{\log n}{\log \log n}$, Proposition~\ref{prop:branching process} implies that with probability $1-o((\log n)^{-3})$, the $h'$-th generation of this process has at least $n^{\gamma}\ge n^{1/2}$ vertices. 
    Hence, using Markov's inequality, with probability $1-o((\log n)^3)$, the number of vertices in $Z_h$ for which less than $n^{1/2}$ vertices in $Z'_{h'}$ satisfy the desired property is at most $|Z_h|/2\le |Z_h|-n^{2/3}/2$, as desired.
    \endproof

    Therefore there exist at least $n^{7/6}/2$ pairs $u^*v^*$ with $v^* \in Z_h$ and $u^* \in Z_{h'}'$ which satisfies the conclusion of Claim~\ref{cl:good}.
    Since these vertices belong to the last generations of the processes, the edges in $R_{i,j}$ between them have not been exposed and, with probability $1-(1-p/K)^{n^{7/6}/2}=1-o((\log n)^{-2})$, one of these pairs appear in $R_{i,j}$.
    This allows us to update $\hat M_i$ by performing the unique sequence of acceptable rotations from $v$ to $v^*$ in $(Z_{\ell})_{\ell\ge 0}$, then the unique sequence of acceptable rotations from $u$ to $u^*$ in $(Z''_{\ell}[v^*])_{\ell\ge 0}$, and finally adding the bridging edge $u^* v^*$. 
    As the first rotation from $v$ was forced to select a pivot strictly outside $C$, this sequence of rotations merges $C$ with at least one other cycle in $M_i \cup M_j$, while at the same time preserving property \ref{item:Q4} since we only used acceptable rotations.
    Thus, $t_s \le t_{s-1} - 1$ and $|W| \le (4(s - 1) + 2)n^{2/3}(\log n)^2 + 2n^{2/3}(\log n)^2 = 4s n^{2/3}(\log n)^2$.

    As $s\le (\log n)^2$, by a union bound, whp our iterative procedure successfully continues until $t_s=1$. 
    It remains to verify that properties \ref{item:Q1}--\ref{item:Q5} are maintained. (Our argument will also show that \ref{item:Q1}, \ref{item:Q2}, \ref{item:Q4} and the stronger version of \ref{item:Q5} are maintained after each iteration, proving \ref{item:Q6} as well.) For \ref{item:Q1} note that, throughout round $(i,j)$, we modify only $M_i^-$. Also, at each iteration we exchange at most $h+h'+1 = O(\log n / \log \log n)$ edges in $\hat M_i$ with edges of $R_{i,j}$ and we make at most $t \le (\log n)^2$ iterations, contributing to a total of at most $O((\log n)^3 / \log \log n) \le (\log n)^3$ exchanges, with each of them contributing at most two edges to the symmetric difference $M_i^+ \Delta M_i^-$, which establishes \ref{item:Q2}.
    Moreover, \ref{item:Q3} is satisfied by design and, since every rotation performed was acceptable and acceptable rotations never split connected components in any 2-factor, \ref{item:Q4} also holds. 
    Finally, we prove the stronger version of \ref{item:Q5} as follows.
    Consider some iteration of the round $(i,j)$ and recall the notation introduced above.
    We let $\hat{M_1}, \dots, \hat{M_k}$ be the matchings at the beginning of the iteration, $e=uv$ be the edge of $\hat{M_i}$ from which we started the branching process, $X$ be the set of used vertices and $\hat{X}$ be the set guaranteed by \ref{item:Q5} and \ref{item:Q6}.
    We have to show that for each $j' \in [j,k]$, each cycle in the union of the $i$-th and $j'$-th matching at the end of the iteration contains a red edge disjoint from $B \cup W$.     
    Recall that the only edge of $\hat{X}$ used in the iteration was $e$ and none of the endpoints of an edge of $\hat{X} \setminus \{e\}$ got added to $W$.
    
    So fix $j' \in [j,k]$ and such a cycle $D$.
    Either $D$ was a cycle at the beginning of the iteration or it was obtained as a result of one or more rotations.
    In the former case, it contains an edge of $\hat{X} \setminus \{e\}$ and this can be used again.
    In the latter case, either $D$ was obtained by merging two or more cycles or the rotations were all performed inside $D$.
    In the first case, consider any two cycles. Then at least one of them contains an edge of $\hat{X} \setminus \{e\}$ and this can be used again.
    In the second case, recall that the first rotation of the iteration was done outside $C$ and $D_{j'}$.
    Therefore if we only performed rotations within $D$, it must mean that $j' >j$ and the cycle containing $e$ in $\hat M_i \cup \hat M_{j'}$ had length at least $n/2$.
    Among its $n/4$ red edges, at most $|B|+2|W| =o(n)$ are edges of $B$ or have endpoints in $W$, and thus there is a suitable choice as well.
    This allows to update the set $\hat{X}$ in all cases and proved \ref{item:Q6}.
    
    Finally, taking a union bound over all $K$ rounds concludes the proof.
\endproof

We now possess the necessary tools to prove Theorem~\ref{thm:main}.
\lateproof{Theorem~\ref{thm:main}}
As mentioned in the introduction, the statement is well-known for $k=1$, so assume $k\ge 2$, $\delta = 1/(8k+2)$, $p=\delta (\log n)/n$ and further fix $\eps,\alpha \in (0,1)$ such that $\eps \ll 1/k$ and $\alpha < \min\{\eps,1/(2k-2)\}$.

By Theorem~\ref{thm:coupling_main}, whp $G(\tau_k)$ can be coupled with $k$ perfect matchings $M_1, \dots, M_k$, a set of vertices~$B$ and a set of edges $R$ such that the properties \ref{item:S1}--\ref{item:S3} hold.
Moreover, the properties \ref{item:P1}--\ref{item:P3} hold as well.

For each $i \in [k]$, define $M_i'$ to be the submatching of $M_i$ consisting of all the edges incident to $B$.
Then $|M_i'| \le n^{1-\eps}$ by property \ref{item:S1} and $M_i'':=M_i \setminus M_i'$ is a uniformly random perfect matching of $V \setminus V(M_i')$, independent of any other colour by property \ref{item:S2}.
Therefore we can apply Lemma~\ref{lem:statistics_new} and thus the family $\{M_1,\dots,M_k\}$ is $\alpha$-shuffled. 
Then we conclude using Lemma~\ref{lemma:rotations}.
\endproof

\section{Concluding remarks}\label{sec:conclusion}

In this paper we established a hitting time result in the random graph process showing that whp~$G(\tau_k)$ contains $k$ edge-disjoint perfect matchings such that any two of them form a Hamilton cycle.
A natural related problem, which to the best of our knowledge has not appeared elsewhere, is the following. Given a perfect matching $M$,
when can one find a second matching $M^*$ such that $M \cup M^*$ is a Hamilton cycle?
Using our tools, one can show that whp~in the random graph process $G(\tau_1)$ contains a perfect matching $M^*$ with the desired property.
It would be interesting to extend this when the number of given matchings is $k$, for any fixed $k$: namely show that given $k$ (edge-disjoint) perfect matchings $M_1,\dots,M_k$, whp~in the random graph process $G(\tau_1)$ contains a perfect matching $M^*$ such that $M_i\cup M^*$ is a Hamilton cycle for every $i\in [k]$.
Our proof mostly works, however, it would be more complicated to establish the statistics we need.
\medskip

Wanless~\cite{wanless:99} formulated an analogue of Kotzig's conjecture in the bipartite setting, stating that $K_{n,n}$ admits a perfect $1$-factorisation if and only if $n=2$ or $n$ is odd.
Analogously, we ask if our results can be adapted to the bipartite setting.
We remark that, if one follows our approach, one would have to apply rotations differently.
Indeed, let $M, M'$ be perfect matchings of $K_{n,n}$ and suppose we would like to merge the cycles in $M \cup M'$ into a Hamilton cycle by applying rotations.
In particular, while building a large set of possible new endpoints, at some point we will have to rotate within a path, say $P$, and say we rotate around its endpoint $v$.
Then we switch some edge $xy$ of $M \cap E(P)$ with a new edge $vx$.
Since we can only use edges which respect the bipartition, $v$ and $x$ must belong to different parts, implying that the rotation will split $P$ into a shorter path and a cycle.
This, at first sight, breaks our strategy as we do not want to increase the number of cycles.
However, it seems plausible that one can still guarantee this by applying more than one rotation at the same time and we leave this as an open direction.
\medskip

Finally, we reiterate the conjecture of Wormald (Conjecture~\ref{conj:wormald}) that for fixed $d \ge 3$, a random $d$-regular graph should whp~have a perfect $1$-factorisation.
We also wonder if it is more generally true that $(n,d,\lambda)$-graphs have a perfect $1$-factorisation (for suitable ranges of parameters).
That they contain a $1$-factorisation has already been proved by  Ferber, Frederickson, Mao, Yepremyan and Zhu~\cite{FFMYZ:26} when $\lambda \le d/12$ and $d$ is large enough.

\section*{Acknowledgments}
Part of this work was done while Lichev was visiting the University of Passau in February 2026, and the author is grateful to the institution for its hospitality.
Glock and Sgueglia are funded by the Deutsche Forschungsgemeinschaft (DFG, German Research Foundation) – 542321564.
During the initial stage of this work, Sgueglia was funded by the Alexander von Humboldt Foundation.
Lichev is supported by the Austrian Science Fund (FWF) [10.55776/ESP624].
For open access purposes, the authors have applied a CC BY public copyright license to any author-accepted manuscript version arising from this submission.

\section*{Statement of AI use}
All mathematical ideas, arguments, results and proofs have been developed and written by the authors.
We used ChatGPT-6 Astra only at the very last stage for local improvements of the~language and style of the manuscript, and for checking the mathematical correctness of all proofs.

\bibliographystyle{amsplain_v2.0customized}
\bibliography{References}

@article{Ros19,
  title={Perfect 1-factorizations},
  author={Rosa, A.},
   JOURNAL = {Math. Slovaca},
  FJOURNAL = {Mathematica Slovaca},
  volume={69},
  number={3},
  pages={479--496},
  year={2019},
  publisher={De Gruyter}
}

@article{Pik19,
  title={A perfect one-factorisation of {$K_{56}$}},
  author={Pike, D. A.},
   JOURNAL = {J. Combin. Des.},
  FJOURNAL = {Journal of Combinatorial Designs},
  volume={27},
  number={6},
  pages={386--390},
  year={2019},
  publisher={Wiley Online Library}
}

@unpublished{GS:25,
    AUTHOR = {Glock, S. and Sgueglia, A.},
     TITLE = {On {K}otzig’s conjecture in random graphs},
   note = {to appear in J. Combin. Theory Ser. B},
}

@incollection {wormald:99,
    AUTHOR = {Wormald, N. C.},
     TITLE = {Models of random regular graphs},
 BOOKTITLE = {Surveys in combinatorics},
    SERIES = {London Math. Soc. Lecture Note Ser.},
    VOLUME = {267},
     PAGES = {239--298},
 PUBLISHER = {Cambridge Univ. Press, Cambridge},
      YEAR = {1999},
}

@article {BF:85,
    AUTHOR = {Bollob\'{a}s, B. and Frieze, A. M.},
     TITLE = {On matchings and {H}amiltonian cycles in random graphs},
     JOURNAL = {Ann. Discrete Math.},
     fjournal={Annals of Discrete Mathematics},
    VOLUME = {28},
     PAGES = {23--46},
      YEAR = {1985},
}

@incollection {kotzig:64,
    AUTHOR = {Kotzig, A.},
     TITLE = {Hamilton graphs and {H}amilton circuits},
 BOOKTITLE = {Theory of {G}raphs and its {A}pplications ({P}roc. {S}ympos.
              {S}molenice, 1963)},
     PAGES = {63--82},
 PUBLISHER = {Publ. House Czech. Acad. Sci., Prague},
      YEAR = {1964},
}

@article{AKS:85,
  title={First occurrence of {H}amilton cycles in random graphs},
  author={Ajtai, M. and Koml{\'o}s, J. and Szemer{\'e}di, E.},
  journal={North-Holland Mathematics Studies},
  volume={115},
  number={C},
  pages={173--178},
  year={1985},
  publisher={Elsevier}
}

@inproceedings {KO:14c,
	AUTHOR = {K\"{u}hn, D. and Osthus, D.},
	TITLE = {Hamilton cycles in graphs and hypergraphs: an extremal
	perspective},
	BOOKTITLE = {Proceedings of the {I}nternational {C}ongress of
	{M}athematicians---{S}eoul 2014. {V}ol. {IV}},
	PAGES = {381--406},
	PUBLISHER = {Kyung Moon Sa, Seoul},
	YEAR = {2014},
}

@incollection {mcdiarmid:89,
    AUTHOR = {McDiarmid, C.},
     TITLE = {On the method of bounded differences},
BOOKTITLE = {Surveys in combinatorics, 1989 ({N}orwich, 1989)},
    SERIES = {London Math. Soc. Lecture Note Ser.},
    VOLUME = {141},
     PAGES = {148--188},
PUBLISHER = {Cambridge Univ. Press},
      YEAR = {1989},
}

@book {JLR:00,
    AUTHOR = {Janson, S. and {\L}uczak, T. and Ruci\'{n}ski, A.},
     TITLE = {Random graphs},
    SERIES = {Wiley-Intersci. Ser. Discrete Math. Optim.},
 PUBLISHER = {Wiley-Interscience},
   ADDRESS = {New York},
      YEAR = {2000},
     PAGES = {xii+333},
}

@article {anderson:73,
    AUTHOR = {Anderson, B. A.},
     TITLE = {Finite topologies and {H}amiltonian paths},
   JOURNAL = {J. Combinatorial Theory Ser. B},
  FJOURNAL = {Journal of Combinatorial Theory. Series B},
    VOLUME = {14},
      YEAR = {1973},
     PAGES = {87--93},
}

@article {wanless:99,
    AUTHOR = {Wanless, I. M.},
     TITLE = {Perfect factorisations of bipartite graphs and {L}atin squares
              without proper subrectangles},
   JOURNAL = {Electron. J. Combin.},
  FJOURNAL = {Electronic Journal of Combinatorics},
    VOLUME = {6},
      YEAR = {1999},
     PAGES = {Research Paper 9, 16},
}

@incollection {bollobas:84,
    AUTHOR = {Bollob\'as, B.},
     TITLE = {The evolution of sparse graphs},
 BOOKTITLE = {Graph theory and combinatorics ({C}ambridge, 1983)},
     PAGES = {35--57},
 PUBLISHER = {Academic Press, London},
      YEAR = {1984},
}

@book{FK:15,
  title={Introduction to random graphs},
  author={Frieze, A. and Karo{\'n}ski, M.},
  year={2015},
  publisher={Cambridge University Press}
}

@article{Luc91,
  title={Size and connectivity of the k-core of a random graph},
  author={{\L}uczak, T.},
   JOURNAL = {Discrete Math.},
  FJOURNAL = {Discrete Mathematics},
  volume={91},
  number={1},
  pages={61--68},
  year={1991},
  publisher={Elsevier}
}

@book{Pen03,
  title={Random geometric graphs},
  author={Penrose, M.},
  volume={5},
  year={2003},
  publisher={OUP Oxford}
}

@article{BBKMW11,
  title={Hamilton cycles in random geometric graphs},
  author={Balogh, J. and Bollob{\'a}s, B. and Krivelevich, M. and M{\"u}ller, T. and Walters, M.},
   JOURNAL = {Ann. Appl. Probab.},
  FJOURNAL = {The Annals of Applied Probability},
  pages={1053--1072},
  year={2011},
  publisher={JSTOR}
}

@article{MPW11,
  title={Disjoint {H}amilton cycles in the random geometric graph},
  author={M{\"u}ller, T. and P{\'e}rez-Gim{\'e}nez, X. and Wormald, N.},
   JOURNAL = {J. Graph Theory},
  FJOURNAL = {Journal of Graph Theory},
  volume={68},
  number={4},
  pages={299--322},
  year={2011},
  publisher={Wiley Online Library}
}

@unpublished{IL25,
  title={Sharp thresholds, hitting times and the power of choice for random geometric graphs},
  author={Ignasiak, D. and Lichev, L.},
  note={arXiv:2512.04191},
  year={2025}
}

@article{Bol90,
  title={Complete matchings in random subgraphs of the cube},
  author={Bollob{\'a}s, B.},
   JOURNAL = {Random Structures Algorithms},
  FJOURNAL = {Random Structures \& Algorithms},
  volume={1},
  number={1},
  pages={95--104},
  year={1990},
  publisher={Wiley Online Library}
}

@article {CE-DGKO24,
    AUTHOR = {Condon, P. and Espuny D\'{\i}az, A. and Gir\~{a}o,
              A. and K\"{u}hn, D. and Osthus, D.},
     TITLE = {Hamiltonicity of random subgraphs of the hypercube},
   JOURNAL = {Mem. Amer. Math. Soc.},
  FJOURNAL = {Memoirs of the American Mathematical Society},
    VOLUME = {304},
      YEAR = {2024},
    NUMBER = {1534},
     PAGES = {v+132},
      ISSN = {0065-9266,1947-6221},
      ISBN = {978-1-4704-7266-5; 978-1-4704-8013-4},
   MRCLASS = {05C45 (05C80)},
}

@article{ER66,
  title={On the existence of a factor of degree one of a connected random graph},
  author={Erd{\H{o}}s, P. and R{\'e}nyi, A.},
  journal={Acta Math. Acad. Sci. Hungar},
  volume={17},
  number={359-368},
  pages={192},
  year={1966}
}

@incollection {BT85,
    AUTHOR = {Bollob\'{a}s, B. and Thomason, A.},
     TITLE = {Random graphs of small order},
 BOOKTITLE = {Random graphs '83 ({P}ozna\'{n}, 1983)},
    SERIES = {North-Holland Math. Stud.},
    VOLUME = {118},
     PAGES = {47--97},
 PUBLISHER = {North-Holland, Amsterdam},
      YEAR = {1985},
      ISBN = {0-444-87821-1},
   MRCLASS = {05C80 (60C05)},
}

@article {FFMYZ:26,
    AUTHOR = {Ferber, A. and Frederickson, B. and Mao, D. and
              Yepremyan, L. and Zhu, Y.},
     TITLE = {Regular bipartite decompositions of pseudorandom graphs},
   JOURNAL = {J. Combin. Theory Ser. B},
  FJOURNAL = {Journal of Combinatorial Theory. Series B},
    VOLUME = {179},
      YEAR = {2026},
     PAGES = {270--292},
      ISSN = {0095-8956,1096-0902},
   MRCLASS = {05C70 (05C80)},
       DOI = {10.1016/j.jctb.2026.04.004},
       URL = {https://doi.org/10.1016/j.jctb.2026.04.004},
}

@unpublished{CS:26+,
    AUTHOR = {Cheng, Y. and Sgueglia, A.},
     TITLE = {The perfect 1-factorisation conjecture holds asymptotically},
   note = {arXiv:2607.09459},
      YEAR = {2026},
}

\end{document}